\documentclass[11pt]{report}

\input xy
\xyoption{all}
\usepackage{epsfig}
\usepackage{graphics}
\usepackage{amsmath}
\usepackage{amscd}
\usepackage{amsfonts,latexsym}

\newenvironment{proof}{\hspace{-\parindent}\textit{Proof.}}{\hfill $\Box$}
\newenvironment{proof2}{\textit{Proof}}{\hfill $\Box$}

\newcommand{\SA}{{\mathcal{A}}}
\newcommand{\SB}{{\mathcal{B}}}
\newcommand{\SC}{{\mathcal{C}}}

\newcommand{\SE}{{\mathcal{E}}}
\newcommand{\SF}{{\mathcal{F}}}

\newcommand{\SI}{{\mathcal{I}}}
\newcommand{\SJ}{{\mathcal{J}}}

\newcommand{\SL}{{\mathcal{L}}}
\newcommand{\SM}{{\mathcal{M}}}

\newcommand{\SO}{{\mathcal{O}}}

\newcommand{\SQ}{{\mathcal{Q}}}

\newcommand{\calS}{{\mathcal{S}}}

\newcommand{\SV}{{\mathcal{V}}}
\newcommand{\SW}{{\mathcal{W}}}

\newcommand{\SZ}{{\mathcal{Z}}}

\newcommand{\PP}{\mathbb{P}}
\newcommand{\ZZ}{\mathbb{Z}}
\newcommand{\CC}{\mathbb{C}}
\newcommand{\QQ}{\mathbb{Q}}

\newcommand{\VV}{\mathbb{V}}
\newcommand{\FM}{\mathfrak{M}}

\newcommand{\OOS}{\SO_S}

\newcommand{\wt}{\widetilde}
\newcommand{\isom}{\cong}
\newcommand{\Ext}{\operatorname{Ext}}
\newcommand{\SExt}{Ext}
\newcommand{\Hom}{\operatorname{Hom}}
\newcommand{\SHom}{Hom}
\newcommand{\Hilb}{\operatorname{Hilb}}
\newcommand{\Quot}{\operatorname{Quot}}
\newcommand{\Sym}{\operatorname{Sym}}
\newcommand{\Pic}{\operatorname{Pic}}
\newcommand{\Spec}{\operatorname{Spec}}

\newcommand{\im}{\operatorname{im}}
\newcommand{\codim}{\operatorname{codim}}
\newcommand{\maxdim}{\operatorname{max.dim}}
\newcommand{\fiber}{\operatorname{fiber}}
\newcommand{\Supp}{\operatorname{Supp}}

\newcommand{\inj}{\hookrightarrow}
\newcommand{\fcndown}[1]{\vcenter{\llap{$\scriptstyle{#1}$}}\Big
\downarrow}
\newcommand{\rk}{\operatorname{rk}}
\newcommand{\gr}{\operatorname{gr}}
\newcommand{\rank}{\operatorname{rank}}

\newcommand{\extlzo}{\Ext ^1 (L \otimes I_Z, \OOS)}
\newcommand{\extn}{\Ext ^1 (L^{\otimes n+1} \otimes I_Z, L^{\otimes  -n})}
\newcommand{\hilbcc}{\Hilb ^{c_2} (S)}
\newcommand{\hilbn}{\Hilb ^{c_2 + n(n+1)L^2} (S)}

\newcommand{\modulil}{\FM (L,c_2)}

\newtheorem{proposition}{Proposition}[section]
\newtheorem{theorem}[proposition]{Theorem}
\newtheorem{definition}[proposition]{Definition}
\newtheorem{lemma}[proposition]{Lemma}

\newtheorem{corollary}[proposition]{Corollary}
\newtheorem{remark}[proposition]{Remark}

\date{23 March 2000
\thanks{Revised version of PhD thesis (Princeton, 1997)}}

\author{Tom{\'a}s L. G{\'o}mez}
\title{Irreducibility of the moduli space of vector bundles on surfaces and
Brill-Noether theory on singular curves
}
\begin{document}

\maketitle

\abstract{
We prove the irreducibility of the moduli space of
rank 2 semistable torsion free sheaves (with a generic polarization and 
any value of $c_2$) on a $K3$ or a del Pezzo  surface. Under these
conditions the moduli space is smooth of the expected dimension.

If the surface $S$ is a $K3$, then by a recent theorem of
O'Grady the moduli space is known to be irreducible. We present a 
new proof of this result using a different technique.

First we consider the case in which $\Pic(S)=\ZZ$. We prove this case
by constructing a connected
family of rank 2 semistable torsion free sheaves that maps onto a
dense set of the moduli space. To prove that this family is connected
we need a result from
Brill-Noether theory. For a smooth curve it is known that the
Brill-Noether locus is connected if the expected dimension is
positive. We need to generalize this for irreducible singular 
curves that lie on a $K3$ surface (we prove it for any surface whose
anticanonical line bundle is generated by global sections).
Finally we remove the condition 
$\Pic(S)=\ZZ$ by considering families of surfaces.

If $S$ is a del Pezzo surface we reduce the problem to the case
of  $\PP^2$ by studying the relationship of moduli spaces
corresponding to different polarizations and then comparing the moduli
space for a surface with the moduli space for the blown up surface at
a point.

}

\tableofcontents

\chapter*{Acknowledgments}

I would like to thank my advisor, Robert Friedman, for suggesting me 
this problem, for many discussions and for his encouragement. I am
very grateful to him for introducing me to the study of moduli spaces
of vector bundles. I consider myself very fortunate for having found
his generous help when I was looking for a thesis problem.

I would also like to thank Ignacio Sols Lucia, who introduced me to 
algebraic geometry. With the time that he
generously dedicated to teach me, he made possible
the transition from my physics undergraduate background to algebraic geometry. 
Thanks for his patience and his friendship.

I have also benefited from discussions with R. Lazarsfeld, R.
MacPherson and many other people as well as from many seminars and
courses at Princeton. 

I want to thank the ``Banco de Espa{\~n}a'' for the fellowship ``Beca
para la ampliaci{\'o}n de estudios en el extranjero'' that supported
my graduate studies. Without their generous support, this thesis
couldn't have been made.

\chapter{Introduction}

In this chapter we will explain the main results of the thesis using
as little mathematical background as possible. We will always work
over the complex numbers, i.e all manifolds will be complex manifolds.
Also we assume that manifolds are projective, i.e. there is an
embedding in $\PP^n_\CC$ (in particular they are K{\"a}hler).

$X$ will be a variety (or a manifold with singularities).
We will consider holomorphic vector bundles over
$X$ (the transition functions are assumed to be holomorphic).
To distinguish non-isomorphic vector bundles of fixed rank 
we can define certain
invariants called Chern classes. These are cohomology classes
$c_i(V)\in H^{2i}(X,\ZZ)$, $1\leq i\leq \dim_\CC X$. We have
$c_i(V)=0$ if $i>\rank(V)$. Even after fixing these discrete
invariants, we can have continuous families of non-isomorphic vector
bundles. I.e., to specify the isomorphism class of a vector bundle it is
not enough to fix some discrete invariants, but we have to fix also
some continuous parameters.  

For fixed rank and Chern classes we would like to define a variety
called the moduli space of vector bundles. Each point of this variety
will correspond to a different vector bundle. In this thesis we will
study the irreducibility of this space for rank 2 vector bundles over
certain surfaces.

Unfortunately, in order to construct this moduli space
we have to restrict our attention to \textit{stable}
vector bundles (this is not a very strong restriction, since it can be
proved that in some
sense all vector bundles can be constructed starting from stable
ones). There are different notions of stability (see chapter
\ref{Preliminaries}). Here we will only discuss \textit{Mumford
stability} (also called slope stability). Let $X$ be a projective
variety and $H$ an ample divisor. For a vector
bundle $V$ we define the slope of $V$ with respect to $H$ as:
$$
\mu_H (V) = \frac { c_1(V) H^{n-1}}{\rank(V)}
$$
(the product is the intersection product or cup product in cohomology)
$V$ is called $H$-stable if for every subbundle $W$ of $V$ we
have
$$
\mu_H (W) < \mu_H (V).
$$
It can be proved that there is a space $\FM^0_H(r,c_i)$, called the
moduli space of $H$-stable vector bundles of rank $r$ and Chern
classes $c_i$ (if the rank is clear from the context we will drop it
from the notation). In many situations it will be a variety (maybe with
singularities) but in general it can have a very singular behavior.

In general $\FM^0_H(r,c_i)$ is not compact. To get a compact moduli
space we need to consider a larger class of objects. Instead of vector
bundles we consider torsion free sheaves (they can be thought as
``singular'' vector bundles that fail to be locally a product on a
subvariety of X). Also we have to relax the stability condition, and
we will consider \textit{Gieseker semistable} sheaves (see chapter
\ref{Preliminaries} for the definition). The moduli space of Gieseker
semistable torsion free sheaves is compact, and we denote it by
$\FM_H(r,c_i)$.

There is an important relationship between the theory of holomorphic
vector bundles and gauge theory: there is a bijection between the set of
$H$-stable vector bundles and differentiable bundles with a
Hermite-Einstein connection. If $c_1=0$ and $\dim_\CC=2$ then a
Hermite-Einstein connection is the same thing as an anti-selfdual (ASD)
connection. Note that the metric of the manifold appears in the ASD
equation. This is reflected in the fact that the stability condition
depends on the polarization. This relationship has been used, for
instance, to calculate Donaldson polynomials for the study of the
topology of four-manifolds.

Now we are going to consider some particular cases. 
Let $C$ be a
curve (i.e., a Riemann surface) of genus $g$. Then the moduli space of
rank $r$ vector
bundles on $C$ is a smooth variety of dimension $(g-1)r^2+1$.

For line bundles (rank=1) over a variety $X$ ($\dim_\CC X=n$)  the
moduli space is called the Jacobian and we have
a very explicit description of it. First we note that
all line bundles are stable (because they don't have subbundle). 
Recall that  $q=h^1(\SO_X)=b_1/2$, where $b_1$ is the first Betti number
(if $X$ is simply connected then $b_1=0$. If $X$ is a curve
then $b_1=2g$). The Jacobian $J$ is the moduli space of line bundles
with $c_1=0$ and it is of the form $\CC^q/\ZZ^{2q}$, where
$\ZZ^{2q}$ is a lattice in $\CC^q$. If $c_1 \neq 0$ then the moduli
space $J^{c_1}$ is isomorphic to the Jacobian, but the isomorphism is
not canonical. 

If $X$ is a curve $c_1$ is called the degree $d$. We define some 
subsets of $J^d$ as follows
$$
W^a_d=\{L \in J^d: \dim(H^0(L))=a+1\}
$$
The study of the properties of these subsets is called Brill-Noether theory.
If the curve is generic, then $W^a_d$ is a subvariety of the
expected dimension $\rho(a,d)=g-(a+1)(g-d+a)$. If $\rho(a,d)>0$ then 
(for any curve) $W^a_d$ is connected. If the curve is
singular, one has to consider also torsion free sheaves in order to get
a compact moduli space (in this case the moduli space will be
singular). This moduli space has been constructed, but little is known
about its Brill-Noether theory.

In this thesis, to
prove the irreducibility of the moduli space of rank 2 vector bundles
on a $K3$ surface, we will need the connectivity of $W^a_d$ for certain 
singular curves that lie
in the surface. In chapter \ref{bn},
theorem I,
we will prove that if $\rho(a,d)>0$, $W^a_d$ is still connected for 
singular curves satisfying certain conditions. 

If $X$ is a complex surface ($\dim_\CC X=2$), then in general the
moduli space of vector bundles of rank $r$ is very singular, but in
many situations it will be a variety (maybe with singularities) of
the expected dimension
$$
\dim
\FM^0_H(r,c_1,c_2)=2rc_2-(r-1)c_1^2-(r^2-1)\chi(\SO_X)+h^1(\SO_X).
$$
By a slight abuse of language we have denoted by $c_2$ the integral 
of the second Chern class
on the variety $\int_X c_2(V)$, and by $c_1^2$ the integral $\int_X
c_1(V) \wedge c_1(V)$. Recall that  $\chi(\SO_X)=1-h^1(\SO_X)+h^2(\SO_X)$.

In this thesis we will consider the case rank=2. For rank 2, 
$\FM^0_H(c_1,c_2)$
is known to be irreducible and of the expected dimension if $c_2>N$,
 where $N$ is a constant that
depends on $X$, $c_1$ and $H$. In chapter \ref{k3} we present a
new proof of the theorem of O'Grady

\smallskip\noindent
\textbf{Theorem II.}
\textit{Let $X$ be a projective $K3$
surface (a $K3$ surface is a simply connected complex surface with
$c_1(T_X)=0$, where $T_X$ is the tangent bundle). Let $H$ be a
generic polarization (see chapter \ref{Preliminaries} for a
definition). Assume that $c_1$ is a nonzero primitive element of $H^2(X,\ZZ)$
(under this condition $\FM^0_H(c_1,c_2)$ is smooth of the
expected dimension).}

\textit{Then $\FM^0_H(c_1,c_2)$ is irreducible.}
\smallskip

We should note that in chapter \ref{k3} we work with the moduli space
of Gieseker semistable torsion free sheaves $\FM_H(c_1,c_2)$.
But it can be shown that under the conditions of the theorem the
points of $\FM_H(c_1,c_2)$ that are not in
$\FM^0_H(c_1,c_2)$ are in a subvariety of positive codimension,
then the irreducibility of one of them is equivalent to the
irreducibility of the other.

Using the relationship of holomorphic vector bundles with gauge theory
this theorem can be stated as follows:

\begin{theorem}
Let $X$ be a projective $K3$ surface with a generic K{\"a}hler
metric $g$. Then the moduli space of
$SO(3)$ anti-selfdual connections (with fixed instanton number $k$ and
second Stiefel-Whitney class $w_2$) is smooth
and connected.
\end{theorem} 

In chapter \ref{dp} we study moduli spaces of vector bundles on
\textit{del Pezzo} surfaces. A del Pezzo surface is a surface whose
anticanonical bundle is ample. It can be shown that these are all the
del Pezzo surfaces: $\PP^2$, $\PP^1\times\PP^1$, or the projective
plane $\PP^2$ blown up at most at 8 generic points. For rank 2 and
fixed Chern classes, the moduli space of stable vector bundles is
known to be empty or irreducible for $X=\PP^2$ or $\PP^1\times\PP^1$.
In chapter \ref{dp} we prove this result for a projective plane
$\PP^2$ blown up at most at 8 generic points, and then we get the
general result:

\smallskip\noindent
\textbf{Theorem III.}
\textit{Let $X$ be a del Pezzo surface with a generic polarization
$H$. Fix some Chern classes $c_1$ and $c_2$. Then
$\FM^0_H(c_1,c_2)$ is irreducible (or empty).}
\smallskip

As in the $K3$ case, in this case irreducibility is equivalent to
connectedness, and we have the same result for $\FM_H(c_1,c_2)$.
We can also translate this theorem into gauge theory language:

\begin{theorem}
Let $X$ be a del Pezzo surface with a generic K{\"a}hler metric $g$.
Then the moduli space of $SO(3)$ anti-selfdual connections (with fixed
Stiefel-Whitney class $w_2$ and instanton number $k$) is
connected (or empty). The same is true for $SU(2)$ anti-selfdual
connections with fixed instanton number $k$.
\end{theorem}

\chapter{Preliminaries}
\label{Preliminaries}

\section{Brill-Noether theory}

Let $C$ be a smooth curve of genus $g$ (we will always assume that the 
base field is $\CC$), $J(C)$ its Jacobian, and $W^r_d(C)$ the
Brill-Noether locus corresponding to line bundles $L$ of degree $d$ and
$h^0(L) \geq r+1$ (see \cite{ACGH}). The expected dimension of this
subvariety is $\rho(r,d)= g -
(r+1)(g-d+r)$. Fulton and Lazarsfeld \cite{F-L} proved that $W^r_d(C)$
is connected
when $\rho > 0$. We are going to generalize this result for certain singular
curves, but before stating our result
(theorem I), we
need to recall some concepts.

Let $C$ be an integral curve (not necessarily smooth). We still have a
generalized Jacobian $J(C)$, defined as
the variety parametrizing line bundles, but it
will not be complete in general. Define the degree of a rank 
one torsion-free sheaf on $C$ to be 
$$
\deg(A)=\chi(A)+p_a-1,
$$
where $p_a$ is the arithmetic genus of $C$.
One can define a scheme $\overline J^d(C)$
parametrizing rank one torsion-free sheaves on $C$ of degree $d$ 
(see \cite{AIK}, \cite {D},  
\cite {R}). If $C$ lies on a surface,
then $\overline
J^d(C)$ is integral, and furthermore the generalized
Jacobian $J(C)$ is an open set in $\overline J^d(C)$, 
and then $\overline J^d(C)$
is a natural compactification of $J(C)$.

We will need to consider
families of sheaves parametrized by a scheme $T$, and furthermore the
curve will vary as we vary the parameter $t\in T$. 

All this can be done
using a relative version of $\overline J^d(C)$, but we will proceed in a
different way. We will use the fact that all these curves are going to
lie on a fixed surface $S$. Then we will think of the coherent sheaves on $C$
as torsion sheaves on $S$ (all sheaves in this paper will be coherent). 
To define precisely which sheaves we will
consider we need some notation. For any sheaf $F$ on $S$, let $d(F)$
be the dimension of its support. We say that $F$ has pure dimension
$n$ if for any subsheaf $E$ of $F$ we have $d(E)=d(F)=n$. Note that if
the support is irreducible, then having pure dimension $n$ is
equivalent to being torsion-free when considered as a sheaf on its
support. The following theorem follows from \cite[theorem 1.21]{S}.

\smallskip
\begin{theorem}[Simpson]
Let $C$ be an integral curve on a surface $S$.
Let $\overline \SJ^d_{|C|}$ be the functor which associates to any scheme
$T$ the set of equivalence classes of sheaves $\SA$ on $S\times T$ with

(a) $\SA$ is flat over $T$.

(b) The induced sheaf $A_t$ on each fiber $S\times \{t\}$ has pure
dimension 1, and its support is an integral curve in the linear system
$|C|$.

(c) If we consider $A_t$ as a sheaf on its support, it is 
torsion-free and has rank one and degree $d$.

Sheaves $\SA$ and $\SB$ are equivalent if there exists a line
bundle $L$ on $T$ such that $\SA \isom \SB \otimes p_T^*L$, where
$p_T:S \times T \to T$ is the projection on the second factor.

Then there is a coarse moduli space that we also denote by $\overline
\SJ^d_{|C|}$. I.e., the points of $\overline \SJ^d_{|C|}$ correspond to
isomorphism classes of sheaves, and for any family $\SA$ of such
sheaves parametrized by $T$, there is a morphism
$$
\phi :T \to \overline \SJ^d_{|C|}
$$
such that $\phi(t)$ corresponds to the isomorphism class of $A_t$.
\end{theorem}

Note that $\overline\SJ^d_{|C|}$ parametrizes pairs $(C',A)$ with $C'$
an integral curve linearly equivalent to $C$ and $A$ a torsion-free
rank one sheaf on $C$.

We denote by $\pi :\overline\SJ^d_{|C|} \to U \subset |C|$ the 
obvious projection giving
the support of each sheaf, where $U$ is the open subset of $|C|$
corresponding to integral curves. 

A family of curves on a surface $S$ parametrized by a curve $T$ is a
subvariety $\SC \subset S \times T$, flat over $T$, such that the
fiber $\SC|_t=C_t$ over each $t\in T$ is a curve on $S$.
Analogously, a family of sheaves on a surface $S$ parametrized by a
curve $T$ is a sheaf $\SA$ on $S\times T$, flat over $T$. For each
$t\in T$ we will denote the corresponding member of the family by
$A_t=\SA|_t.$

Altman, Iarrobino and Kleiman
\cite{AIK} proved the following theorem

\begin{theorem}[Altman--Iarrobino--Kleiman]
With the same notation as before, $\overline \SJ^d_{|C|}$ 
is flat over $U$ and its
geometric fibers are integral. The subset of $\overline \SJ^d_{|C|}$
corresponding to line bundles (i.e., the relative generalized Jacobian) is
open and dense in $\overline \SJ^d_{|C|}$.
\end{theorem}

We also
consider the family of generalized Brill-Nother loci $\overline
\SW^r_{d,|C|} \subset \overline \SJ^d_{|C|}$, and the projection 
$q:\overline\SW^r_{d,|C|} \to U$.

We can define the generalized Brill-Noether locus $\overline W^r_d(C)$
as the set of points in $\overline J^d(C)$ corresponding to sheaves $A$
with $h^0(A) \geq r+1$ (note that it is complete because of the
upper semicontinuity of $h^0(\cdot)$). There is also a determinantal
description that
gives a scheme structure. This description is a straightforward 
generalization of the description for smooth curves (see \cite
{ACGH}), but we are only interested in the connectivity of $\overline
W^r_d(C)$, so we can give it the reduced scheme structure. 

We will consider curves that lie on a surface $S$ with the following
property:

\bigskip
$h^1(\OOS)=0$ \textit{, and }$-K_S$ \textit{ is generated by global
sections.}\hfill (*)
\bigskip

We will need this condition to prove proposition \ref{bn1.5}.
For instance, $S$ can be a K3 surface or a del Pezzo surface with
$K_S^2\neq 1$.
Now we can state the theorem that we are going to prove in chapter 
\ref{bn}.

\smallskip
\noindent\textbf{Theorem I.}
\textit{ Let $C$ be a reduced irreducible curve of arithmetic genus
$p_a$ that lies in a surface $S$
satisfying  \textup{(*)}.
Let $\overline J^d(C)$, $d>0$, be the compactification of the generalized
Jacobian. Then for any $r \geq 0$ such that $\rho(r,d)= p_a -
(r+1)(p_a-d+r)>0$, the generalized Brill-Noether subvariety $\overline
W^r_d(C)$ is nonempty and connected.}
\smallskip

\section{Moduli space of torsion free sheaves}

Let $X$ be a smooth projective variety over $\CC$ of dimension $n$. We are
interested in
constructing a moduli space that will parametrize torsion free sheaves
on $X$, with some fixed rank $r$ and Chern classes $c_i$. To do so, we
have to introduce some notion of stability. 

\begin{definition}\textup{\textbf{(Mumford-Takemoto stability)}}
Fix a polarization $H$. For any nonzero torsion free sheaf $F$ we define
the slope of $F$ with respect to $H$ as
$$
\mu_H(F)={\frac {c_1(F)\cdot H^{n-1}} {\rk (F)}}
$$ 
We will say that a torsion free sheaf $V$ on $X$ is Mumford H-stable
(resp. semistable) if for every nonzero subsheaf $W$ of $V$ we have
$$
\mu_H(W) < \mu_H(V) \ \ \ (\text{resp. }\leq ).
$$
\end{definition}

Using this definition one can construct the coarse moduli space of
Mumford stable vector bundles. This notion of stability turns out to
have an analog in differential geometry: Mumford stable holomorphic
vector bundles are in one to one correspondence with differentiable
vector bundles with a Hermite-Einstein connection (the choice of a
polarization is replaced by the choice of a Riemannian metric). This
correspondence
was proved by Narasimhan and Seshadri \cite{N-S} if $X$ is a curve,
by Donaldson \cite{Do1,D-K} for a surface, and it was then generalized for
any dimension by \cite{U-Y}. This correspondence has been useful to
calculate differentiable invariants of 4-manifolds (the so called
Donaldson invariants).

The moduli space of Mumford stable sheaves is in general not compact. To
define a compactification we have to introduce a refined notion of
stability.

\begin{definition}\textup{\textbf{(Gieseker stability).}}
Fix a polarization $H$. For any nonzero torsion free sheaf $F$ define
the Hilbert polynomial of $F$ with respect to $H$ as
$$
p_H(F)(n) =\frac {\chi(F \otimes \SO_X(nH))}{rank(F)}.
$$
Given two polynomials $f$ and $g$, we will write $f\prec g$ (resp.
$\preceq$) if $f(n)<g(n)\ ($resp. $\leq)$ for $n\gg0$.

We will say that $V$ is Gieseker stable (resp. semistable) if for
every nonzero torsion free subsheaf $W$ we have
$$
p_H(W) \preceq p_H(V) \ \ \ (\text{resp. }\preceq ).
$$
\end{definition}

Using Hirzebruch-Riemann-Roch theorem we can see that Gieseker
semistability implies Mumford semistability, and Mumford stability
implies Gieseker stability.

In order to have a separated moduli space it is not enough to consider
isomorphism classes of sheaves. We will introduce the notion of
S-equivalence. For any Gieseker semistable sheaf $V$ there is a
filtration (\cite{G,Ma})
$$
0 =V_0 \subseteq V_1\subseteq \cdots \subseteq V_t=V
$$
such that $V_i/V_{i-1}$ is stable and $p_H(V_i)=p_H(V)$. We define
$\gr (V)=\oplus (V_i/V_{i-1})$. It can be proved that $\gr(V)$ doesn't 
depend on
the filtration chosen. We will say that $V$ is S-equivalent to $V'$ if
$\gr(V)=\gr(V')$.

There is another characterization of S-equivalence that is more
illuminating from the point of view of moduli problems. Assume that we
have a family of Gieseker semistable sheaves parametrized by a curve
$T$. I.e., we have a sheaf $\SV$ on $X \times T$, flat over $T$
inducing torsion free
Gieseker stable sheaves $\SV|_t$ on the slices $X \times \{t\}$. Assume
that for one point $0\in T$ we have $\SV|_0\isom V$ and for the rest of
the points $\SV|_t$ is isomorphic to some other fixed $V'$. We will
say that $V$ and $V'$ are equivalent. The equivalence relation
generated by this definition is S-equivalence.

It can be proved that there is a coarse moduli space for
S-equivalence classes of Gieseker semistable torsion free sheaves
(with fixed rank and Chern classes). This moduli space is projective.

In general the moduli space can be very singular, but if $X$ is a surface
and for fixed rank $r$, $c_1$ and polarization, it is
known that for $c_2$ large enough the singular locus is a proper
subset of positive codimension of the moduli space
\cite{Do2,F,Z,G-L2}. The moduli space has the expected dimension
$$
2rc_2-(r-1)c_1^2-(r^2-1)\chi(\SO_X),
$$
and is irreducible
\cite{G-L1,G-L2,O1,O2}. In the rank two case it is also known, again
for $c_2$ large enough, that the moduli space is normal and has local
complete intersection singularities at points corresponding to stable
sheaves, and if the surface $X$ is of general type (with some
technical conditions), then also the moduli space is of general type
\cite{L2}.

It is natural to ask what is the effect of changing the choice of
polarization. From now on we will assume that $X$ is a surface $S$ and
that the rank is 2 unless otherwise stated. We will denote the moduli
space of rank 2 torsion free sheaves that are Gieseker semistable with
respect to the polarization $H$ by $\FM_H(c_1,c_2)$.

\begin{definition}
Fix $S$, the first and second Chern classes $c_1$, $c_2$.
Let $\zeta$ be some class in $H^2(S,\ZZ)$ with
$$
\zeta\equiv c_1\ (\text{mod}\ 2),\ \ c_1^2-4c_2\leq \zeta^2 <0.
$$
The wall of type $(c_1,c_2)$ associated to $\zeta$ is a hyperplane of
$H^2(S,\QQ)$ with
nonempty intersection with the ample cone $\Omega_S$
$$
W^\zeta=\{x\in \Omega_S |\  x \cdot \zeta =0\}. 
$$
The connected components of the complement of the walls in the ample
cone are called chambers.
\smallskip
A polarization $H$ is called $(c_1,c_2)$-generic if it doesn't lie on
a wall (i.e., it lies in a chamber).
\end{definition}

The walls of type $(c_1,c_2)$ are known to be locally finite on the
ample cone \cite{F-M}.

If a polarization is $(c_1,c_2)$-generic it is easy to see that
Mumford and Gieseker stability coincide, and furthermore there are no
strictly semistable sheaves. Stable sheaves are \textit{simple}
($\Hom(V,V)=\CC$). If $-K_S$ is effective then this fact and the
Kuranishi local model for the moduli space proves that the moduli space is
smooth of the expected dimension (if not empty) \cite{F}. 

If $H_1$ and $H_2$ are two generic polarizations in the same chamber,
then every sheaf that is $H_1$-stable is also $H_2$-stable, and we can
identify the corresponding moduli spaces \cite{F,Q1,Q2}.

If we restrict our attention to some particular class of algebraic
surfaces we can obtain more properties of the moduli space (without
the condition on $c_2$). 

The moduli space of sheaves on a $K3$ surface with a generic polarization
has been studied by Mukai
\cite{M} when the expected dimension is 0 or 2. In particular he
proved that the moduli space is irreducible. O'Grady \cite{O3} has proved
irreducibility for any $c_2$ (and any rank), as well as having
obtained results about the Hodge structure.
In chapter \ref{k3} (theorem II),  we give a new
proof of the irreducibility for any $c_2$ (and rank 2) based on our
results about Brill-Noether theory on singular curves.

Now let's consider the case $S=\PP^2$ and rank 2.
Tensoring with a line bundle we can assume that $c_1$ is either $0$ or
$1$. If $c_1=0$, then the moduli space is empty for $c_2<2$ and is
irreducible of the expected dimension for $c_2\geq2$. If $c_1=1$ then
the moduli space is empty for $c_2<1$ and irreducible of the expected
dimension for $c_2\geq1$.

In the case $S=\PP^1\times\PP^1$ (and also rank 2), if we take a
$(c_1,c_2)$-generic polarization, the moduli space is also known to be
either empty or irreducible. We will generalize this result for any
del Pezzo surface (i.e., a surface with $-K_S$ ample) in chapter 
\ref{dp} (theorem III).

\chapter{Connectivity of Brill-Noether loci for singular curves}
\label{bn}

\setcounter{proposition}{0}

Recall (see chapter \ref{Preliminaries}) that we are going to study the
Brill-Noether locus of singular irreducible curves that lie on a
smooth surface $S$ satisfying

\bigskip
$h^1(\OOS)=0$ \textit{, and }$-K_S$ \textit{ is generated by global
sections.}\hfill (*)
\bigskip

Now we state the theorem that we are going to prove:

\smallskip
\noindent\textbf{Theorem I.}
\textit{ Let $C$ be a reduced irreducible curve of arithmetic genus
$p_a$ that lies in a surface $S$
satisfying  \textup{(*)}.
Let $\overline J^d(C)$, $d>0$, be the compactification of the generalized
Jacobian. Then for any $r \geq 0$ such that $\rho(r,d)= p_a -
(r+1)(p_a-d+r)>0$, the generalized Brill-Noether subvariety $\overline
W^r_d(C)$ is nonempty and connected.}
\smallskip

\begin{remark}
\label{remark}
\textup{If $r \leq d-p_a$, by Riemann-Roch inequality we have 
$\overline W^r_d(C) 
= \overline
J^d(C)$, and this is connected. Then, in order to prove theorem I we 
can assume
$r>d-p_a$. Note that if $A$ corresponds to a point in $\overline
W^r_d(C)$ with $r>d-p_a$, then by Riemann-Roch theorem $h^1(A) > 0$.}
\end{remark}

\bigskip
\textbf{Outline of the proof of theorem I}
\bigskip

Note that $\overline W^r_d(C)$ is the fiber of $q$ over the point $u_0
\in |C|$ corresponding to the curve $C$.
Let $U$ be the open subset of $|C|$ corresponding to integral curves,
and $V$ the subset of smooth curves. Define $(\overline
\SW^r_d)_V$ to be the Brill-Noether locus of sheaves with smooth
support, i.e. $(\overline\SW^r_d)_V=q^{-1}(V)$. 
 By \cite{F-L}, the restriction
$q^{}_V: (\overline \SW^r_d)_V \to V$ has connected fibers. We want to
use this fact to show that $\overline W^r_d(C)$ is connected. Let $A$ 
be a rank one
torsion-free sheaf on $C$ corresponding to a point in $\overline
W^r_d(C)$,
and assume that it is generated by global sections. We think of $A$
as a torsion sheaf on $S$. Then we have a
short exact sequence on $S$
$$
0 \to E \stackrel{f_0} \to H^0(A) \otimes \OOS \to A \to 0,
$$ 
where the map on the right is evaluation. This sequence has already
appeared in the literature (see \cite{La}, \cite{Ye}).
Our idea is to deform $f_0$ to a family $f_t$. The cokernel of $f_t$ 
will define
a family of sheaves $A_t$ with $h^0(A_t) \geq h^0(A)$ (because
$h^0(E)=0$), and then for each $t$ the point in $\overline
\SJ^d_{|C|}$ corresponding to $A_t$ lies in $\overline \SW^r_{d,|C|}$. 
Assume that there are 'enough' homomorphisms from $E$ to $H^0 \otimes
\OOS$ and the family $f_t$ can be chosen
general enough, so that for a general $t$, the support of $A_t$ is
smooth (the details of this construction are in section
\ref{bnParticular case}). The family $A_t$ shows that the point
in $\overline
W^r_d(C)$ corresponding to $A$ is in the closure of $(\overline
\SW^r_d)_V$ in $\overline \SJ^d_{|C|}$. It can be shown that this 
closure has connected
fibers. Let $X$ be the fiber over $u_0$ of this closure. Then all sheaves 
for which this construction works are in the
connected component $X$ of $\overline
W^r_d(C)$. If this could be
done for all sheaves in $\overline
W^r_d(C)$ this would finish the proof, but there are sheaves for which
this construction doesn't work. For these sheaves we show in
section \ref{bnGeneral case} that they can be deformed (keeping the
support $C$ unchanged) to a sheaf for which a refinement of this 
construction works.
This shows that all points in $\overline W^r_d(C)$ are in the connected
component $X$.

\section{The main lemma}
\label{bnMain lemma}

The precise statement that we will use to prove theorem I is
the following lemma.

\begin{lemma}
\label{bn0.2}

Let $C$ be an integral complete curve in a surface $S$. Assume that for each 
rank one
torsion-free sheaf $A$ on $C$ with $h^0(A)=r+1>0$ and $\deg(A)=d>0$
such that $\rho(r,d)>0$ we
have the following data:

\smallskip

(a) A family of curves $\SC$ in $S$ parametrized by an irreducible
curve $T$ (not necessarily complete).

(b) A connected curve $T'$ (not necessarily irreducible nor complete) 
with a map $\psi:T' \to T$.

(c) A rank one torsion-free sheaf $\SA$ on $\SC' = \SC \times _T T'$,
flat over $T'$,
inducing rank one torsion-free sheaves on the fibers of $\SC' \to T'$.

\smallskip

Assume that the following is satisfied:

\smallskip

(i) $\SC|_{t^{}_0} \isom C$ for some $t^{}_0 \in T$, $\SC|_t$ is linearly 
equivalent to $C$ for all $t\in T$, and $\SC|_t$
is smooth for $t \neq t^{}_0$.

(ii) One irreducible component of $T'$ is a finite cover of $T$, 
and the rest of the components of $T'$ are mapped to $t^{}_0 \in T$.

(iii) $\SA|_{t'_0} \isom A$ for some $t'_0 \in T'$ mapping to 
$t^{}_0 \in T$.

(iv) $h^0(\SA|_{t'}) \geq r+1$ for all $t' \in T'$.

\smallskip

Then the generalized
Brill-Noether subvariety $\overline W^r_d(C)$ of the compactified
generalized Jacobian $\overline J^d(C)$ is connected.

\end{lemma}

\begin{proof}
We will use the notation introduced in the previous section. The map
$q:\overline \SW^r_{d,|C|} \to U$ is a projective morphism. Recall
that $\overline W^r_d(C)$ is the fiber of $q$ over $u_0$, where $u_0$ 
is the point corresponding to $C$. By
\cite{F-L} the morphism $q$ has connected
fibers over $V$, thus a general fiber of $q$ is connected, and we
want to prove that the fiber over $u_0 \in U$ is also
connected.
$$
\begin{array}{ccc}
\overline W^r_d(C) & \inj & \overline \SW^r_{d,|C|} \\
\fcndown{}         &  {}  &      \fcndown{q}         \\
u_0                & \inj &       U                 
\end{array}
$$
Let $\overline \SW^r_{d,|C|} \stackrel{q'} \to U' \stackrel{g} \to U$ be
the Stein factorization of $q$ (see \cite[III Corollary 11.5]{H}),
 i.e. $q'$ has connected fibers and
$g$ is a finite morphism. A general fiber of $q$ is connected, and
then $U'$ has one irreducible component $Z$ that maps to $U$
birationally.
The subset $U$ is open in $|C|$ and hence normal, the restriction 
$g|_Z:Z \to U$ 
is finite and birational, $Z$ and
$U$ are integral, thus by
Zariski's main theorem (see \cite[III Corollary 11.4]{H}) each fiber
of $g|_Z$ consists of just one point. Let
$z_0$ be the point of $Z$ in the fiber $g^{-1}(u_0)$.

\textbf{\textit{Claim.}} Let $y_0$ be a point in the fiber
$q^{-1}(u_0)=
\overline W^r_d(C)$. Then $y_0$ is
mapped by $q'$ to $z_0$.

This claim implies that that $\overline
W^r_d(C)$ is connected. Now we will prove the claim.

Let $A$ be the sheaf on $S$ corresponding to the point $y_0$. Let $T'$, $T$,
$t'_0\in T'$, $t^{}_0\in T$, $\psi:T' \to T$ be the
curves points and morphism given by the hypothesis of the lemma. 
Let $\phi:T' \to
\overline\SJ^d_{|C|}$ be the morphism given by the universal property
of the moduli space $\overline\SJ^d_{|C|}$. Item (iv) imply that the
image of $\phi$ is in $\overline \SW^r_{d,|C|}$. 

\centerline{
\xymatrix{
{}\save[]+<0.7cm,0cm>*{t'_0\in} \restore
& T' \ar[d]_\psi \ar[r]^\phi & \overline\SW^r_{d,|C|}
 \ar[d]^{q'} \ar@(dl,ul)[dd]_q
& {}\save[]+<-0.7cm,0cm>*{\ni y_0} \restore\\
{}\save[]+<0.7cm,0cm>*{t^{}_0\in} \restore& T & U' \ar[d]^{g}& 
Z \ar@{_{(}->}[l] \ar[dl]^{g|_Z} & {}\save[]-<0.6cm,0cm>*{\ni z^{}_0} \restore\\
& & U & {}\save[]+<-0.8cm,0cm>*{\ni u^{}_0} \restore}
}

\begin{figure}[ht]

\centerline{\epsfig{file=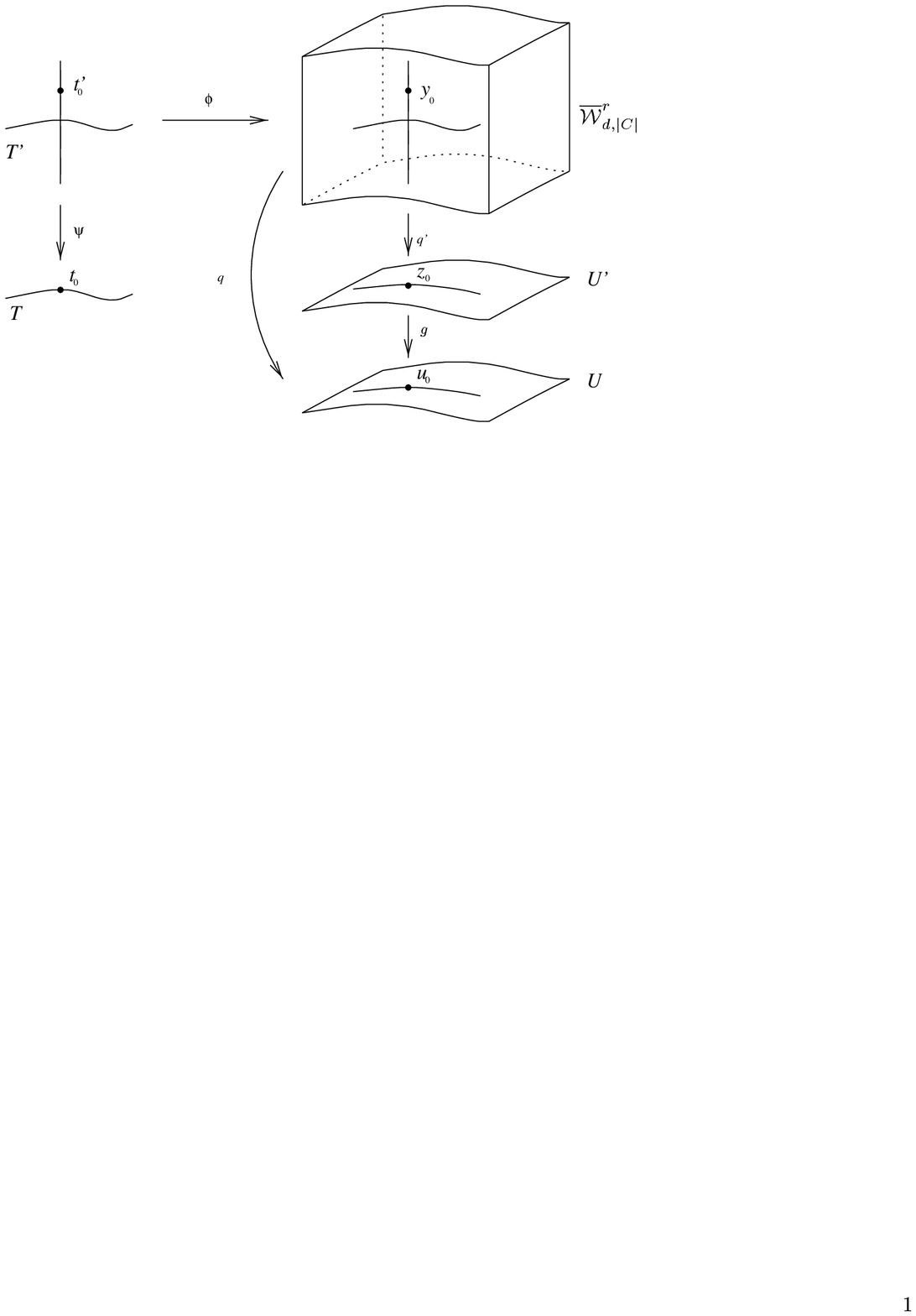,height=3in,width=4.5in,
bbllx=0in,bblly=9.187in,bburx=4.5in,bbury=12.187in,clip=}}
\end{figure}

The restriction of $q' \circ \phi$ to $T' \setminus \psi^{-1}(t^{}_0)$
maps to $Z$, because for $t' \in T' \setminus
\psi^{-1}(t^{}_0)$ the sheaf $\SA|_{t'}$ has smooth support by item (i).
Items (c) and (i) imply that $g\circ q' \circ \phi
(\psi^{-1}(t^{}_0))=u_0$ . Thus $q' \circ \phi (\psi^{-1}(t^{}_0))$ is a
finite number of points (because it is in the fiber of $g$ over $u_0$).

The facts that $q' \circ \phi(T' \setminus \psi^{-1}(t^{}_0))$ is in
$Z$ and that $q' \circ \phi (\psi^{-1}(t^{}_0))$ is a
finite number of points imply that $q' \circ \phi (\psi^{-1}(t^{}_0))$ is
also in $Z$ (because by item (b) the curve $T'$ is connected and thus also its
image under $q'\circ \phi$), and in fact $q' \circ \phi
(\psi^{-1}(t^{}_0))=z_0$ because $q' \circ \phi (\psi^{-1}(t^{}_0))$ is in
the fiber of $g$ over $u_0$.

By item (ii), $t'_0 \in \psi^{-1}(t^{}_0)$. Then $q' \circ \phi
(t'_0)=z_0$, and by item (iii)  we have $y_0=\phi(t'_0)$, then
$q'(y_0)= q'(\phi(t'_0))=z_0$ and the claim is proved.

\end{proof}

In section \ref{bnParticular case} we will construct this family under 
some assumptions on $A$ (proposition \ref{bn1.5}),
and in section \ref{bnGeneral case} we will show how to use that to 
construct a family for any $A$. Note that because of remark \ref{remark} we can
assume  $h^1(A)>0$.

\section{A particular case}
\label{bnParticular case}

Given a rank one torsion-free sheaf $A$ on an integral curve lying on
a surface $S$, we
define another sheaf $A^*$ that is going to be some sort of dual. Let
$j$ be the inclusion of the curve $C$ in the surface $S$. We
define $A^*$ as follows:
$$ 
A^*= \SExt^1(j_*A,\omega_S).
$$
The operation $A \to A^*$  is a contravariant functor. 
Note that the support of $A^*$ is
$C$. It will be clear from the
context when we are referring to $A^*$ as a torsion sheaf on $S$ or as a sheaf
on $C$. In the case in which $A$ is a line bundle, then $A^*=A^\vee
\otimes \omega_C$. Now we prove some properties of this ``dual''.

\begin{lemma}
\label{bn1.1}
Let $A$ be a rank one torsion-free sheaf on an integral curve lying on
a surface. Then $A^{**}=A$
\end{lemma}

\begin{proof}
First observe that if $L$ is a line bundle on $C$, then $(A \otimes
L)^* \isom A^* \otimes L^\vee$. 
To see this, take an injective resolution of $\omega_S$
$$
0 \to \omega_S \to \SI_0 \to \SI_1 \to \cdots
$$
Now we use this resolution to calculate the $\SExt$ sheaf.
$$
\SExt^1(A\otimes L,\omega_S)  =  h^1(\SHom(A \otimes L,\SI_\bullet))=
h^1(L^\vee\otimes \SHom(A,\SI_\bullet)) =
$$
$$
 =L^\vee\otimes h^1
(\SHom(A,\SI_\bullet))=L^\vee\otimes \SExt^1(A,\omega_S)
$$
The third equality follows from the fact that $\SHom(A,\SI_\bullet)$
is supported on the curve and $L^\vee$ is locally free

It follows that $(L\otimes A)^{**}\isom L\otimes A^{**}$, and then
proving the lemma for $A$ is equivalent to proving it for $L\otimes A$.
Multiplying with an appropriate very ample line bundle, we can assume
that $A$ is generated by global sections. Then we have an exact
sequence
\begin{equation}
0 \to E \to V \otimes \OOS \to A \to 0,
\label{eqbn1.1}
\end{equation}
where $V=H^0(A)$. The following lemma proves that $E$ is locally free.

\begin{lemma}
\label{bn1.2}
Let $M$ be a torsion-free sheaf on an integral curve $C$ that lies on a
smooth surface $S$. Let $j:C \to S$ be the inclusion. Let $F$ be a
locally free sheaf on the surface. Let
$f:F \to j_* M$ be a surjection. Then
the elementary transformation $F'$ of $F$, defined as the kernel of
$f$
\begin{equation}
\label{eqbn1.1bis}
0 \to F' \to F \stackrel{f}{\to} j_* M \to 0,
\end{equation}
is a locally free sheaf.
\end{lemma}

\begin{proof}
$M$ is torsion-free sheaf on $C$, and then $j_* M$ has depth at least one,
and because $S$ is smooth of dimension 2, this implies that the
projective dimension of $j_* M$ is at most one 
($\SExt^i(j_* M,\SO_S)=0$
for $i \geq 2$). Now $\SExt^i(F,\SO_S)=0$ for $i \geq 1$ because $F$ is
locally free, and then from the
exact sequence \ref{eqbn1.1bis}, we get
$$
0 \to \SExt^i(F',\OOS) \to \SExt^{i+1}(j_* M,\OOS) \to 0,
\quad i\geq 1,
$$ 
and then $\SExt^i(F',\SO_S)=0$ for $i \geq 1$, and
this implies that $F'$ is locally free.

\end{proof}

In particular, $E^{\vee\vee}=E$. Applying the
functor $\SHom(\cdot,\omega_S)$ twice to the sequence \ref{eqbn1.1}, we get
$$
0 \to E \to V \otimes \OOS \to A^{**} \to 0.
$$ 
Comparing with \ref{eqbn1.1} we get the result (because the map on the
left is the same for both sequences).

\end{proof}

\begin{lemma}
\label{bn1.3}
$\Ext^1(A,\omega_S) \isom H^0(A^*)$, and this is dual to $H^1(A)$.
\end{lemma}

\begin{proof}
The local to global spectral sequence for Ext gives the following
exact sequence
$$ 
0 \to H^1(\SHom(A,\omega_S)) \to \Ext^1(A,\omega_S) \to H^0(A^*)
\to H^2(\SHom(A,\omega_S))
$$
But $\SHom(A,\omega_S)=0$ because $A$ is supported in $C$ and then the first
and last terms in the sequence are zero and we have the desired
isomorphism.

\end{proof}

Now we will prove a lemma that we will need. The proof can also be
found in \cite{O}, but for convenience we reproduce it here.

\begin{lemma}
\label{bn1.4}
Let $E$ and $F$ be two vector bundles of rank $e$ and $f$ over a
smooth variety $X$. Assume that
$E^\vee \otimes F$ is generated by global sections. 
If $\phi: E \to F$ is a sheaf morphism, we define $D_k(\phi)$ to be
the subset of X where $\rk(\phi_x) \leq k$ (there is an obvious
determinantal description of $D_k(\phi)$ that gives a scheme
structure). Let $d_k$ be the expected dimension of $D_k(\phi)$
$$
d_k=\dim (X) - (e-k)(f-k).
$$
Then there is a Zariski dense set $U$ of $\Hom(E,F)$ such that if
$\phi\in U$, then 
we have that $D_k(\phi) \setminus
D_{k-1}(\phi)$ is smooth of the expected dimension (if $d_k < 0$ then it
will be empty).
\end{lemma}

\begin{proof}
Let $M_k$ be the set of matrices of dimension $e \times f$ and of rank at
most $k$
(there is an obvious determinantal description that gives a
scheme structure to this subvariety). It is well known that the
codimension of $M_k$ in the space of all matrices is
$(e-k)(f-k)$, and that the singular locus of $M_k$ is $M_{k-1}$.

Now, because $E^\vee \otimes F$ is generated by global sections,
we have a surjective morphism
$$
H^0(E^\vee \otimes F) \otimes \SO_X \to E^\vee \otimes F
$$
that gives a morphism of maximal rank between the varieties defined as
the total space of the previous vector bundles
$$
p: X \times H^0(E^\vee \otimes F) \to
\VV (E^\vee \otimes F).
$$
Define $\Sigma_k \subset \operatorname{V}(E^\vee \otimes F)$ as
the set such that $\rk(\phi_x) \leq k$. The fiber of $\Sigma_k$ over
any point in $X$ is obviously $M_k$. Define $Z_k$ to be $p^{-1}(\Sigma_k)$. The
fact that $p$ has maximal rank implies that $Z_k$ has codimension
$(e-k)(f-k)$ in $X \times H^0(E^\vee \otimes F)$ and that the
singular locus of $Z_k$ is $Z_{k-1}$.

Now observe that the restriction of the projection
$$
q|_{Z_k \setminus Z_{k-1}} :Z_k \setminus Z_{k-1} \to H^0(E^\vee
\otimes F)
$$
has fiber ${q|_{Z_k \setminus Z_{k-1}}}^{-1}(\phi) \isom D_k(\phi) \setminus
D_{k-1}(\phi)$. Finally, by generic smoothness, for a general $\phi \in 
H^0(E^\vee \otimes F)$ this is smooth of the expected dimension
(or empty).

\end{proof}

Now we will construct the deformation of $A$ that we described in
the section \ref{bnMain lemma}
in the particular case in which both $A$ and $A^*$ are generated by
global sections. 

\begin{proposition}
\label{bn1.5}
Let $A$ be a rank one torsion-free sheaf on an integral curve $C$
lying on a surface $S$ with $h^1(\OOS)=0$ and $-K_S$ generated by
global sections. Denote $j:C \inj S$. If $A$ and $A^*$ are both generated by
global sections, then there exists a
(not necessarily complete) smooth irreducible curve $T$ and a sheaf
$\SA$ on $S
\times T$ flat over $T$, such that

\smallskip
(a) the sheaf
induced on the fiber of $S \times T \to T$ over some $t^{}_0 \in T$ is $j_* A$

(b) the sheaf $A_t$ induced on the fiber over any $t\in T$ with $t \neq 
t^{}_0$ is supported on a
smooth curve $C_t$ and it is a rank one torsion-free sheaf when 
considered as a sheaf on $C_t$

(c) $h^0(A_t) \geq h^0(A)$ for every $t\in T$.
\end{proposition}

Note that these are the hypothesis of lemma \ref{bn0.2} for the particular
case in which both $A$ and $A^*$ are generated by global sections. We
will lift this condition in the next section.

\begin{proof}
The fact that $A$ is generated by global sections implies that there is
an exact sequence
$$
0 \to E \stackrel{f_0}{\to} V \otimes \OOS \to A \to 0 \; \; \; \; \; 
V=H^0(A),
$$
with $E$ locally free (by proposition \ref{bn1.2}). Taking global sections in
this sequence we see that $H^0(E)=0$, because
$$
0 \to H^0(E) \to V \stackrel{\isom}{\to} H^0(A).
$$
Consider a curve $T$ mapping to $\Hom(E,V \otimes \OOS)$ with $t^{}_0 \in
T$ mapping to $f_0$ (so that item (a) is satisfied). Denote by $f_t$ 
the morphism given for $t \in T$
by this map. After shrinking $T$ we can assume that $f_t$ is still
injective. Let $\pi_1$ be the projection of $S \times T$ onto the
first factor and let $\SE=\pi_1^*E$. Using the universal sheaf and
morphism on $\Hom(E,V \otimes
\OOS)$ we can construct (by pulling back to $S \times T$) an exact
sequence on $S \times T$
$$
0 \to \SE \stackrel{f}{\to} V \otimes \SO_{S \times T} \to \SA \to 0
$$
that restricts for each $t$ to an exact sequence
\begin{equation}
\label{eqbn1.2}
0 \to E \stackrel{f_t}{\to} V \otimes \OOS \to A_t \to 0,
\end{equation}
where $A_t$ is a sheaf supported in the degeneracy locus of $f_t$. It
is clear that $\deg(A)=\deg(A_t)$.

Now we are going to prove that if the curve $T$ and the mapping to
$\Hom(E, V \otimes \OOS)$ are chosen generically, the quotient of the map
gives the desired deformation.

The flatness of $\SA$ over $T$ follows from the fact that it has a
short resolution and from the local criterion of flatness (we
can apply \cite[III Lemma 10.3.A]{H}).  

The condition on $h^0(A_t)$ follows because $H^0(E)=0$ and we have a
sequence
$$
0 \to H^0(E)=0 \to V \to H^0(A_t),
$$ 
and then $h^0(A) \leq h^0(A_t)$. This proves item (c).

Using the long exact sequence obtained
by applying $\SHom(\cdot,\OOS)$ to \ref{eqbn1.2}, and the fact that $E$ is
locally free, we obtain that
$\SExt^i(A_t,\OOS)$ vanishes for $i \geq 2$, and so the projective dimension of
$A_t$ is 1, and this implies that $A_t$, when considered as a sheaf on
its support $C_t$, is torsion-free.

We have to prove that we can choose the curve $T$ and the map
to $\Hom(E,V \otimes \OOS)$ such that $C_t$ is smooth for $t \neq t^{}_0$
(here we will use that $A^*$ is
generated by global sections).

First note that $\SExt^1(A,\OOS)$ is generated by global sections,
because $\SExt^1(A,\OOS)=A^*\otimes \omega_S^{-1}$, and both $A^*$ and
$\omega_S^{-1}$ are generated by global sections.
Now we see that $E^\vee$ is generated by global sections, because
we have
$$
0 \to V^\vee \otimes \OOS \to E^\vee \to \SExt^1(A,\OOS) \to 0,
$$
$\SExt^1(A,\OOS)$ is generated by global sections and $H^1(V^\vee \otimes
\OOS)=0$. Then $E^\vee \otimes (V \otimes \OOS)$ is generated by
global sections.

Now apply lemma \ref{bn1.4} with $F=V \otimes \OOS$. Then $n=m=r+1$, 
$k=r$ and the
expected dimension is 1. And the lemma gives that for $\phi$ in a
Zariski open subset of
$\Hom(E, V\otimes \OOS)$, the degeneracy locus $D_r(\phi)$ of $\phi$ is
smooth away from the locus $D_{r-1}(\phi)$ where $\phi$ has rank
$r-1$, but again by
lemma \ref{bn1.4} the locus $D_{r-1}(\phi)$ is empty. This proves item
(b).

\end{proof}

\section{General case}
\label{bnGeneral case}

Now we don't assume that $A$ satisfies the properties of the
particular case (i.e., $A$ and $A^*$ now might not be generated by
global sections). We will find a new sheaf that satisfies those
conditions. We know how to deform this new sheaf, and we will show how
we can use this deformation to construct a deformation of the original $A$.

We start with a rank one torsion-free sheaf $A$ with $h^0(A)$,
$h^1(A)>0$ on an integral curve
$C$ lying on a surface.
First we define $A'$ as the base point free part of $A$, i.e. $A'$ is
the image of the evaluation map
$$
H^0(A)\otimes \SO_C \to A.
$$
We have assumed that $h^0(A)>0$, and then $A'$ is a (nonzero) rank one
torsion-free sheaf.
Obviously, $H^0(A)=H^0(A')$. We have a short exact sequence
$$
0 \to A' \to A \to Q \to 0,
$$
where $Q$ has support of dimension 0. Now consider ${A'}^*$, 
and define $B$ to be
its base point free part. We have $h^0({A'}^*)=h^1(A')=h^1(A)+h^0(Q)
\geq h^1(A)>0$. The first equality by lemma \ref{bn1.3}, and the last
inequality by assumption. Then $B$ is a (nonzero) rank one 
torsion-free sheaf. Finally define $A''$ to be equal to $B^*$.

\begin{lemma}
\label{bn2.1}
Both $A''$ and ${A''}^*$ are generated by global sections.
\end{lemma}

\begin{proof}
Since $B$ is the base point free part of ${A'}^*$, we have a sequence
$$
0 \to B \to {A'}^* \to R \to 0
$$
where $R$ has support of dimension zero. Applying $\SHom(\cdot,\omega_S)$
we get
$$
0 \to A' \to B^*=A'' \to \wt R \to 0 \; \; \; \; \; \wt
R=\SExt^2(R,\omega_S),
$$
whose associated cohomology long exact sequence gives
$$
0 \to H^0(A') \to H^0(B^*) \to H^0(\widetilde R) \to H^1(A') \to
H^1(B^*) \to 0.
$$
To see that $A''$ is generated by global sections, it is enough to
prove that the last map is an isomorphism, because
then the first three terms make a short exact sequence, and the fact
that $A'$ and $\widetilde R$ are generated by global sections (the first
by definition, the second because its support has dimension zero) will
imply that $B^*$ (that is equal to $A''$ by definition) is generated
by global sections.

To prove that the last map is an isomorphism, we only need to show
that $h^1(A')=h^1(B^*)$, and this is true because
$$
h^1(A')=h^0({A'}^*)=h^0(B)=h^0(B^{**})=h^1(B^*).
$$
The first equality is by lemma
\ref{bn1.3}, the second because $B$ is the base point free part of 
${A'}^*$, the
third by lemma \ref{bn1.1},  and 
the last again by lemma \ref{bn1.3}.

To see that ${A''}^*$ is generated by global sections, note that by
definition ${A''}^*=B^{**}=B$, and this is generated by global sections.

\end{proof}

We started with a rank one torsion-free sheaf $A$ with $h^0(A)$ and
$h^1(A)>0$, and we
have constructed new sheaves $A'$ and $A''$ with (nontrivial) maps $A'
\to A$ and $A' \to A''$. They give rise to exact sequences
\begin{eqnarray}
0 \to A' \to A \to Q \to 0 \nonumber \\
\label{eqbn2.1}
0 \to A' \to A'' \to \widetilde Q \to 0
\end{eqnarray}

\begin{lemma}
\label{equalit}
With the previous definitions we have $h^0(A')=h^0(A)$ 
and $h^1(A'')=h^1(A')$.
\end{lemma}

\begin{proof}
By construction $h^0(A')=h^0(A)$ and $h^0({A''}^*)=h^0({A'}^*)$. By
lemma \ref{bn1.3} this last equality is equivalent to
$h^1(A'')=h^1(A')$.

\end{proof}

As $A''$ and ${A''}^*$ are generated by global sections, then by
proposition \ref{bn1.5} the sheaf $A''$ can be
deformed in a family $A''_t$ in such a way that the support of a 
general member of the
deformation is smooth. The idea now is to find (flat) deformations of
$A'$ and $A$, so that for every $t$ we still have maps like \ref{eqbn2.1}.
{}From the existence of these maps we will be able to obtain the
condition that $h^0(A_t) \geq h^0(A)$,
then we will be able to apply lemma \ref{bn0.2} and then
theorem I will be proved. The details are in 
section \ref{bnProof of the main theorem}. We will start by showing 
how the condition on
$h^0(A_t)$ is obtained, and then how we can find the deformations of
$A'$ and $A$.

\begin{proposition} 
\label{bn2.2}
Let $A$, $A'$, $A''$ be rank one torsion-free sheaves on an integral
curve $C$. Assume that they fit into exact sequences like \ref{eqbn2.1}
and that $h^0(A')=h^0(A)$ and $h^1(A')=h^1(A'')$.
Let $P$ be a curve (not necessarily complete), and
let $\SA$, $\SA'$, and $\SA''$ be sheaves on $S \times P$, flat over
$P$, inducing for each
$p \in P$ rank one torsion-free sheaves $A^{}_p$, $A'_p$, $A''_p$,
supported on a curve $C_p$ of $S$, where $A^{}_{p^{}_0}=A$, $A'_{p^{}_0}=A'$, 
and
$A''_{p^{}_0}=A''$ for some $p^{}_0 \in P$. Assume that $h^0(A''_p) \geq 
h^0(A''_{p^{}_0})$ for all $p\in P$ and that we have
short exact sequences 
$$
0 \to \SA' \to \SA \to \SQ \to 0
$$
$$
0 \to \SA' \to \SA'' \to \widetilde \SQ \to 0
$$
with $\SQ$ and $\widetilde \SQ$ flat over $P$ (i.e., the induced
sheaves $Q_p$, $\widetilde Q_p$ have constant length, equal to $l(Q)$
and $l(\widetilde Q)$ respectively).

Then we have $h^0(A^{}_p) \geq h^0(A^{}_{p^{}_0})$ for all $p\in P$.
\end{proposition}

\begin{proof}
For each $p \in P$ we have sequences
$$
0 \to A'_p \to A^{}_p \to Q^{}_p \to 0
$$
$$
0 \to A'_p \to A''_p \to \widetilde Q^{}_p \to 0.
$$
The maps on the left are injective because they are nonzero and
the sheaves have rank one and are torsion-free.
Using the associated long exact sequences and the hypothesis we
have
$$
h^0(A^{}_p) \geq h^0(A'_p) \geq h^0(A''_p)-l(\widetilde Q^{}_p) \geq h^0(A'')- 
l(\widetilde Q)=h^0(A')=h^0(A).
$$

\end{proof}

It only remains to prove that those sheaves can be ``deformed along'',
and that those deformations are flat, i.e. that given $A$, $A'$ and
$A''$ we can construct $\SA'$ and $\SA''$. This is proved in the following
propositions.

\begin{proposition}
\label{bn2.3}
Let $L$ and $M$ be rank one torsion-free sheaves on an integral curve
$C$ that lies on a surface $S$. Assume we have
a short exact sequence

\begin{equation}
\label{eqbn2.2}
0 \to L \to M \to Q \to 0.
\end{equation}
Assume furthermore that we are given a sheaf $\SM$ on $S \times P$
(where $P$
is a connected but not necessarily irreducible curve) that is a
deformation of $M$, flat over $P$. I.e., $\SM|_{p^{}_0}
\isom M$ for some $p^{}_0 \in P$, and for all $p\in P$ we have that
$M_p=\SM|_p$ are torsion-free 
sheaves on $C_p$,
where $C_p$ is a curve on $S$.

Then, there is a connected curve $P'$ 
with a map $f:P' \to P$ and a
sheaf $\SL'$ over $S \times P'$ with the following properties:

One irreducible component of $P'$ is a finite cover of $P$ and the
rest of the components map to $p^{}_0 \in
P$. The sheaf $\SL'$ is a deformation of $L$, in the sense that 
$\SL'|_{p'_0} \isom L$
for some $p'_0 \in P'$ mapping to $p^{}_0 \in P$, the sheaf $\SL'$ is flat over
$P'$ and induces rank one torsion-free sheaves on the fibers over
$P'$. And if we define $\SM'$ to be the pullback of $\SM$ to $S
\times P'$, there exists an exact sequence
$$
0 \to \SL' \to \SM' \to \SQ' \to 0,
$$
inducing short exact sequences
$$ 
0 \to L'_{p'} \to M'_{p'} \to Q'_{p'} \to 0
$$
for every $p'\in P'$.
\end{proposition}

\begin{proof}
If the support of $Q$ were in the
smooth part of the curve, we would have $M \isom L \otimes \SO_C(D)$, with
$D$ an effective divisor of degree $l(Q)$. Then, if we are given a
deformation $M_p$ of $M$, we only need to find a deformation $D_p$ of the
effective divisor $D$, with the only condition that $D_p$ is an
effective divisor on $C_p$, with degree $l(Q)$. This can easily be
done if we are in the analytic category. In general we might need to
do a base change of the parametrizing curve $P$ and we will obtain a 
finite cover $P'$ of $P$ (What we are doing is
moving a dimension zero and length $l(Q)$ subscheme of $S$, with the
only restriction that for each $p$ the corresponding scheme is in $C_p$). 
Then we only need to define $L_{p'}=M_{p'} \otimes
\SO_{C_{p'}}(-D_{p'})$ and the
proposition would be proved (with $P'$ a finite cover of $P$).

To be able to apply this, we will have to make first a deformation of
$L$, keeping $M$ fixed, until we get $Q$ to be supported in the smooth
part of $C$ (the curve $C$ also remains fixed in this deformation).
This is the reason for the need of the curve $P'$ with some
irreducible components mapping to $p^{}_0$.

We will prove this by induction on the length of the intersection of
the support of $Q$ and the singular part of $C$. 

\begin{lemma}
\label{bn2.4}
Let $L$ and $M$ be rank one torsion-free sheaves on an integral curve
$C$ that lies on a surface $S$. Assume we have
a short exact sequence
$$
0 \to L \to M \to Q \to 0.
$$
Assume that $Q=R \oplus Q'$ where $Q'$ has length $l(Q)-1$ and it is
supported in the smooth part of $C$, and $R$ has length one as it is
supported in a singular point of $C$ (``the length of the 
intersection of the
support of $Q$ and the singular part of $C$ is one'').

Then there is a flat deformation $L_y$ of $L$ parametrized by a
connected curve $Y$ (it might not be 
irreducible) such that $L_{y_0}=L$ for
some $y_0\in Y$ and for every $y\in Y$ there is an exact sequence
$$
0 \to L_y \to M \to Q_y \to 0
$$
and there is some $y_1\in Y$ such that the support of $Q_{y_1}$ is in
the smooth part of $C$.
\end{lemma}

\begin{proof}
In this situation, the exact sequence \ref{eqbn2.2}
gives rise to another exact sequence
$$
0 \to L \otimes I_Z^\vee \to M \to R \to 0
$$
where the map on the right is the composition of $M \to Q$ and the
projection $Q \to R$, and we denote by $I_Z$ the ideal sheaf of the
support $Z$ of $Q'$.
Because $Z$ is in the smooth part of $C$, $I_Z$ is an invertible sheaf.
Note that $Q'$ is the quotient of $\SO_C$ by this ideal sheaf.
Define $\widehat L$ to be $L \otimes I_Z^\vee$. If we know how to
make a flat deformation $\widehat L_y$ of $\widehat L$ so that the quotient
$R_y$ is
supported in the smooth part of $C$ for some $y_1\in Y$,
then we can construct a deformation $L_y$ of $L$ defined as
$$
L_y = \widehat L_y \otimes I_Z.
$$
Note that this deformation is also flat. The cokernel $Q_y$ of $L_y
\to M$ is supported in the smooth part of $C$ for the
points $y\in Y$ for which $R_y$ is supported in the smooth part of $C$.

This shows that to prove the lemma we can assume that $Q$ has length
one and its support is a singular point of $C$, i.e. $Q=\SO_x$, where
$x$ is a singular point of $C$.

Consider the scheme $\Quot^1(M)$ representing the functor of quotients
of $M$ of length 1. If the support $x$ of the quotient $Q$ is in the
smooth part of $C$, then there is only one surjective map (up to
scalar) because $\dim \Hom (M,Q)=1$, whose kernel is $M\otimes
\SO_C(-x)$.

If $x$ is in the singular part, then in general $\dim\Hom(M,Q)>1$, and
the quotients are parametrized by $\PP\Hom(M,Q)$ (the universal bundle
is flat
over $\PP\Hom(M,Q)$). We
want to show that $\Quot^1(M)$ is connected by constructing a flat
family of quotients $M \to \wt Q_{\tilde c}$ (the family $\wt
Q_{\tilde c}$ will 
be parametrized by an open set of the normalization $\wt C$ of $C$)
such that for a general ${\tilde c}$ the support
of $\wt Q_{\tilde c}$ is in the smooth part of $C$, and for some point 
$\tilde c_0$ the support of
$\wt Q_{\tilde c_0}$ is a singular point of $C$.

Consider the normalization $\widetilde C$ of $C$, and let $F$ be an
open set of $\wt C$
\begin{equation}
\label{eqbn2.3}
\CD
F @>j>> C \times F \\
@.              @V{\pi_1}VV      \\
@.                 C             \\
\endCD
\end{equation}
Where $\pi_1$ is the projection to the first factor and $j=(\nu,i)$,
the morphism
$\nu:F\inj \wt C\to C $ being the restriction to $F$ of the
normalization map and $i$ the
identity map. Note that $j$
is a closed immersion, and its image is just $C \times_C F \isom F$.

Let $\tilde c_0$ be a point of $\widetilde C$ in $\nu^{-1}(x)$ (the
family is
going to be parametrized by an open neighborhood $F$ of $\tilde c_0$).
We have to
construct a surjection of $\wt \SM=\pi_1^* M$ onto $\wt \SQ = j_*
\SO_{F}$. Note that $\wt\SQ|_{C \times \tilde c }=\wt
Q_{\tilde c} \isom \SO_{\nu
(\tilde c)}$ and that
$\wt\SQ$ is flat over $F$.

Now, to define that quotient, it is enough to define it in the
restriction to the image of $j$ (because this is exactly
the support of $\wt \SQ$). So the map we have to define is
$$
j^*\wt\SM \to \SO_{F} \; .
$$
But $j^*\wt\SM=\nu^*M$ is a rank one sheaf on the smooth curve
$F$, so it is the direct sum of a line bundle and a torsion
part $T$. Shrinking $F$ if necessary, the line bundle part is 
isomorphic to $\SO_F$, and we have
$$ 
j^*\wt\SM \isom T \oplus \SO_F \; , 
$$
and then to define the quotient we just take an isomorphism in the
torsion-free part. This finishes the proof of the lemma.

\end{proof}

Now we go to the general case: the intersection of the support of $Q$
with the singular part of $C$ has length $n$. We are going to see how
this can be reduced to the case $n=1$.

Take a surjection from $Q$ to a sheaf $Q'$ of length $n-1$, such that
$Q$ is isomorphic to $Q'$ at the smooth points. The kernel $R$ of this
surjection will have length 1, and will be supported in a singular
point of $C$. It is isomorphic to $\SO_x$, for some singular point $x$. 
We have a diagram
$$
\CD
@.      @.     0   @.  0      @.    \\
@.      @.    @VVV   @VVV     @.    \\
0  @>>> L @>>> L'@>>>  R @>>> 0     \\
@.      @|    @VVV   @VVV     @.    \\
0  @>>> L @>>> M @>>>  Q @>>> 0     \\
@.      @.    @VVV   @VVV     @.    \\
@.      @.     Q' @=   Q'     @.    \\
@.      @.    @VVV   @VVV     @.    \\
@.      @.     0   @.  0      @.    \\
\endCD
$$
Observe that $L$, $L'$ and $R$ satisfy the hypothesis of lemma \ref{bn2.4},
 so we
can find deformations $L_y$, $R_y$ (parametrized by some curve $Y$ and
with $L_{y_0}=L$ and $R_{y_0}=R$ for some $y_0 \in Y$) such 
that for some $y_1 \in Y$ we have that the support of the 
corresponding sheaf $R_{y_1}$ is
a smooth point of $C$. All the maps of the previous diagram can be
deformed along. To do this, we change $L$ by $L_y$, $R$ will be
deformed to $R_y$ and $L'$ is kept constant. Then $Q$ is deformed to a
family $Q_y$ defined as $M/L_y$. The cokernel of $R_y \to Q_y$ will be
$Q_y/R_y=M/L'=Q'$, and hence we keep it constant.
Then for each $y$ we still have a commutative
diagram, and furthermore it is easy to see that all deformations
are flat (note that $R_y$ is a flat deformation and $Q'$ is kept
constant, and then $Q_y$ is a flat deformation). An important point is
that $M$ remains fixed, and the
injection $L \to M$ is deformed to $L_y \to M$.
$$
\CD
@.      @.     0   @.  0      @.    \\
@.      @.    @VVV   @VVV     @.    \\
0  @>>> L_y @>>> L'@>>>  R_y @>>> 0     \\
@.      @|    @VVV   @VVV     @.    \\
0  @>>> L_y @>>> M @>>>  Q_y @>>> 0     \\
@.      @.    @VVV   @VVV     @.    \\
@.      @.     Q' @=   Q'     @.    \\
@.      @.    @VVV   @VVV     @.    \\
@.      @.     0   @.  0      @.    \\
\endCD
$$
For $y_1$ we have that the length of the intersection
of the support
of $Q_{y_1}$ with the singular part of $C$ is $n-1$. We repeat the process
(starting now with $L_{y_1}$, $M$ and $Q_{y_1}$), until all the points of the
support of $Q$ are moved to the smooth part of $C$. This finishes
the proof of the proposition.

\end{proof}

The following proposition is similar to proposition \ref{bn2.3}, 
but now the roles of
$L$ and $M$ are changed: we are given a deformation of $L$ and we have
to deform $M$ along.

\begin{proposition}
\label{bn2.5}

Let $L$ and $M$ be rank one torsion-free sheaves on an integral curve
$C$ that lies on a surface $S$. Assume we have
a short exact sequence
\begin{equation}
\label{eqbn2.4}
 0 \to L \to M \to Q \to 0.
\end{equation}
Assume furthermore that we are given a sheaf $\SL$ on $S \times P$
(where $P$ is a connected but not necessarily irreducible curve) that 
is a deformation of $L$, flat over $P$ , i.e., $\SL|_{p^{}_0}
\isom L$ for some $p^{}_0\in P$, and for all $p\in P$, we have that 
$L_p=\SL|_p$ are torsion-free sheaves on $C_p$,
where $C_p$ is a curve on $S$.
 
Then, there is a connected curve $P'$
with a map $f:P' \to P$ and a
sheaf $\SM'$ over $S \times P'$ with the following properties:

One irreducible component of $P'$ is a finite cover of $P$ and
the rest of the components map to $p^{}_0 \in
P$. The sheaf $\SM'$ is a deformation of $M$, in the sense 
that $\SM'|_{p'_0} \isom M$
for some $p'_0 \in P'$ mapping to $p^{}_0 \in P$, the sheaf $\SM'$ is
flat over
$P'$ and induces rank one torsion-free sheaves on the fibers over
$P'$. And if we define $\SL'$ to be the pullback of $\SL$ to $S
\times P'$, there exists an exact sequence
$$
0 \to \SL' \to \SM' \to \SQ' \to 0,
$$
inducing short exact sequences
$$ 
0 \to L'_{p'} \to M'_{p'} \to Q'_{p'} \to 0
$$
for every $p'\in P'$.

\end{proposition}

\begin{proof}
The proof is very similar to the
proof of proposition \ref{bn2.3}. Again we start by observing that if the
support of $Q$
were in the smooth part of the curve, we would have $M \isom L \otimes
\SO_C (D)$, with $D$ an effective divisor. Then if we are given a
flat deformation $L_p$ of $L$, we find a deformation $D_p$ of $D$ as
in the first part, and the proposition would be proved. So again we
need a lemma that deforms $Q$ so that its support is in the smooth
part of $C$.

\begin{lemma}
\label{bn2.6}
Let $L$ and $M$ be rank one torsion-free sheaves on an integral curve
$C$ that lies on a surface $S$. Assume we have
a short exact sequence
$$
0 \to L \to M \to Q \to 0.
$$
Assume that the part of $Q$ with support in the smooth part of $C$ has
length $l(Q)-1$, i.e. $Q=R \oplus Q'$, where $R$ has length one and is
supported in a singular point of $C$ and $Q'$ has length $l(Q)-1$ and
is supported in the smooth part of $C$. Then there is a flat
deformation $M_y$ of $M$ parametrized by a curve $Y$, such that for 
every $y\in Y$ there is an exact sequence
$$
0 \to L \to M_y \to Q_y \to 0
$$
with $M_y$ a torsion-free sheaf, and there is some $y_1 \in Y$ such
that the support of $Q_{y_1}$ is in
the smooth part of $C$.
\end{lemma}

\begin{proof}
Arguing as in the proof of lemma \ref{bn2.4}, we see that it is enough to
prove the case $l(Q)=1$, and $Q=\SO_x$ for $x$ a singular point of
$C$, then we can assume that the extension of the hypothesis of the
lemma is 
\begin{eqnarray}
0 \to L \to M \to \SO_x \to 0.
\label{assumeext}
\end{eqnarray}
Now we will consider all extensions of $\SO_x$ (for $x$ any point in $C$)
by $L$. If $x$ is a smooth point, then there is only one extension
that is not trivial (up to equivalence)
$$
0 \to L \to M \isom L\otimes \SO_C(x) \to \SO_x \to 0.
$$
All these extensions are then parametrized by the smooth part of $C$.

But if $x$ is a singular point, we could have more extensions, because
in general
$s=\dim \Ext^1(\SO_x,L)>1$. They will be parametrized by a projective
space $\PP ^{s-1}$. We call this space $E_x$.
Note that there is a universal extension on $C \times E_x$ that is
flat over $E_x$.
 We denote by $e_1$ the
point in $E_x$ corresponding to the extension \ref{assumeext}.

Assume that $\wt Q_{\tilde c}$ is a
family of torsion sheaves on $C$ with length 1, 
parametrized by a curve $F$ such that for a general point
$\tilde c\in F$ of the parametrizing curve the support of $\wt Q_{\tilde
c}$ in $C$ is a smooth point, and
for a special point $\tilde c_0\in F$ the support of
$\wt Q_{\tilde c_0}$ is a
singular point.
Now assume that we can construct a flat family (parametrized by $F$) of
nontrivial 
extensions of $\wt Q_{\tilde c}$ by $L$. The extension corresponding
to $\tilde c_0$ gives a point $e_2$ in $E_x$. The space $E_x$ is a
projective space, thus connected, and then there is a curve containing
$e_1$ and $e_2$. Using this curve (together with the universal
extension for $E_x$) and the curve $F$ (together with the family of
extensions that it parametrizes) we construct the
curve $Y$ that proves the lemma.

Now we need to construct $F$. As in the proof of lemma \ref{bn2.4}, 
the parametrizing curve $F$ will be an affine neighborhood of $\tilde
c_0$ in the normalization $\wt C$ of $C$, where $\tilde c_0$ is a
point that maps to the singular point $x$ of $C$. Consider again the
diagram \ref{eqbn2.3} of the proof of lemma \ref{bn2.4}. The
family will be given by an extension of $\wt\SQ =
j_*\SO_F$ by $\wt\SL=\pi^*_1 L$ on $C \times F$. These extensions are
parametrized by the group $\Ext^1(\wt\SQ,\wt\SL)$. The following lemma gives
information about this group and relates this extension with the
extensions that we get after restriction for each slice $C \times
\tilde c$.
We will call $\wt Q_{\tilde c}$ and $\wt L_{\tilde c}$ the restrictions of
$\wt\SQ$ and $\wt\SL$ to
the slice $C \times {\tilde c}$. Note that the restriction $\wt
L_{\tilde c}$ is isomorphic to $L$.

\begin{lemma}
\label{bn2.7}
With the previous notation, we have 

1) $\Ext^1(\wt{\SQ},\wt\SL) \isom H^0(\SExt^1(\wt\SQ,\wt\SL))$

2) $\SExt^1(\wt\SQ,\wt\SL)$ has rank zero outside of the support of
$\wt\SQ$, and
rank 1 on the smooth points of the support of $\wt\SQ$

3) Let $I$ be the ideal sheaf corresponding to a slice $C \times
{\tilde c}$. Then the natural map 
$$
\SExt^1_{\SO_{C \times F}} (\wt\SQ,\wt\SL) \otimes \SO_{C \times F}/I \to
\SExt^1_{\SO_{C \times {\tilde c}}}(\wt Q_{\tilde c},\wt L_{\tilde c})
$$
is injective.

\end{lemma}

\begin{proof}
Item 1 follows from the fact that
$\SHom(\wt\SQ,\wt\SL)=0$ and the exact sequence
$$
0 \to H^1(\SHom(\wt\SQ,\wt\SL)) \to \Ext^1(\wt\SQ,\wt\SL) \to
H^0(\SExt^1(\wt\SQ,\wt\SL)) \to H^2(\SHom(\wt\SQ,\wt\SL)).
$$

To prove item 2 note
that the stalk of $\SExt^1(\wt\SQ,\wt\SL)$ at a point $p$ is isomorphic to 
$\Ext^1(R/I,R)$, where $R$ is the local ring
at the point $p$, and $I$ is the ideal defining the support or $\wt\SQ$.
The ideal $I$ is principal if the point $p$ is smooth, then $R/I$ has
a free resolution
$$
0 \to I \to R \to R/I \to 0
$$
and it follows that $\Ext^1(R/I,R)\isom R/I$.

For item 3, consider the exact sequence 
$$
0 \to \wt\SQ \stackrel{\cdot f}{\to} \wt\SQ \to \wt Q_{\tilde c} \to 0
$$
where the first map is multiplication by the local equation $f$ of the
slice $C \times {\tilde c}$. Applying $\SHom(\cdot,\wt L_{\tilde c})$ we get
$$
\SHom(\wt\SQ,\wt L_{\tilde c})=0 \to \SExt^1(\wt Q_{\tilde c},\wt
L_{\tilde c}) \to 
\SExt^1(\wt\SQ,\wt L_{\tilde c}) \to
\SExt^1(\wt\SQ,\wt L_{\tilde c}),
$$
but the last map is zero. To see this, take a locally free resolution
of $\wt\SQ$. The map induced on the resolution by the multiplication with
the equation $f$ is just multiplication by the same $f$ on each term 
$$
\CD
\SF^\bullet @>>> \wt\SQ @>>> 0          \\
@V{\cdot f}VV @V{\cdot f}VV @.   \\
\SF^\bullet @>>> \wt\SQ @>>> 0          \\
\endCD
$$
A local section of the sheaf $\SExt^i(\wt\SQ,\wt L_{\tilde c})$ is
represented by some
local section $\varphi(\cdot)$ of $\SHom(\SF^i,\wt L_{\tilde c})$, 
and the endomorphism
induced by multiplication by $f$ on $\SExt^i(\wt \SQ,\wt L_{\tilde
c})$ is given by
precomposition with multiplication $\varphi(f\cdot)$, but $\varphi$ is
a morphism of sheaves of modules and then this is equal to $f
\varphi(\cdot)$, and this is equal to zero because $f \wt L_{\tilde
c}=0$.
Then we have that 
\begin{equation}
\label{iso}
\SExt^1(\wt Q_{\tilde c},\wt L_{\tilde c}) \isom 
\SExt^1(\wt\SQ,\wt L_{\tilde c}).
\end{equation}
Taking the exact sequence
$$
0 \to \wt\SL \stackrel{\cdot f}{\to} \wt\SL \to \wt L_{\tilde c} \to 0
$$
and applying $\SHom(\wt\SQ,\cdot)$ we get
$$
\SExt^1(\wt\SQ,\wt\SL) \stackrel{\cdot f}{\to} \SExt^1(\wt\SQ,\wt\SL) \to 
\SExt^1(\wt\SQ,\wt L_{\tilde c})
$$
and using this and the isomorphism \ref{iso} we have an injection
$$
\SExt^1(\wt\SQ,\wt\SL) \otimes \SO_{C \times F}/I \isom 
\SExt^1(\wt\SQ,\wt\SL)/(f\cdot \SExt^1(\wt\SQ,\wt\SL)) \inj
\SExt^1(\wt Q_{\tilde c},\wt L_{\tilde c}).
$$

\end{proof}

Now we are going to construct the family of extensions. By item 2 of
the lemma 
the sheaf $\SE=\SExt^1(\wt \SQ,\wt \SL)$ is isomorphic to
$\SO_X \oplus T(\SE)$ (shrinking $F$ if necessary) where $X$ is the 
support of $\wt\SQ$ and $T(\SE)$
is the torsion part. Take a nonvanishing section of the torsion-free
part, and by item 1 this gives a nonzero element $\psi$ of
$\Ext^1(\wt\SQ,\wt\SL)$. This element gives a nontrivial extension
$$
0 \to \wt\SL \to \wt\SM \to \wt\SQ \to 0.
$$
Observe that $\wt\SM$ is flat over the base, because both $\wt\SL$ and
$\wt\SQ$ are
flat. 

By items 3 and 1 we have that the image of $\psi$ under the
restriction map 
$$
\Ext^1(\wt\SQ,\wt\SL) \to \Ext^1(\wt Q_{\tilde c},L)
$$
is nonzero for any ${\tilde c}$ (recall that $\wt L_{\tilde c}=L$ for all 
${\tilde c}$), and this means that the
extensions that we
obtain after restriction to the corresponding slices 
\begin{equation}
\label{eqbn2.5}
0 \to L \to \wt M_{\tilde c} \to \wt Q_{\tilde c} \to 0
\end{equation}
are non trivial. Furthermore $\wt M_{\tilde c}$ is torsion-free. To prove this
claim, let $T(\wt M_{\tilde c})$
be the torsion part of $\wt M_{\tilde c}$. The map $L \to T(\wt
M_{\tilde c})$ coming
from \ref{eqbn2.5} is zero, because $L$ is torsion-free, i.e.
$T(\wt M_{\tilde c})$ injects in $\wt Q_{\tilde c}$. Then we have
$$
\wt Q_{\tilde c} \isom {\frac {\wt M_{\tilde c}} L} \isom 
{\frac{\wt M_{\tilde c}/T(\wt M_{\tilde c}) \oplus T(\wt M_{\tilde c})} L}
\isom {\frac{\wt M_{\tilde c}/T(\wt M_{\tilde c})} L} \oplus T(\wt
M_{\tilde c}).
$$
$\wt Q_{\tilde c}$ doesn't decompose as the direct sum of two sheaves,
and then one of these summands must be zero.
The first summand cannot be zero, because this would
imply that $L \isom \wt M_{\tilde c}/T(\wt M_{\tilde c})$ and then 
$\wt M_{\tilde c} \isom L \oplus
\wt Q_{\tilde c}$, contradicting the hypothesis that the extension is
not trivial.
Then we must have $T(\wt M_{\tilde c})=0$, and the claim is proved.

\end{proof}

Now we are going to consider the general case, in which the part of
$Q$ supported in singular points has length $n$. We are going to see
that this can be reduced to the case $n=1$, in a similar way to
proposition \ref{bn2.3}.

Let $R=\SO_x$, where $x$ is a singular point in the support of $Q$,
and take a surjection from $Q$ to $R$. We have a diagram
$$
\CD
@.      @.     @.     0     @.     \\
@.      @.     @.    @VVV   @.     \\
@.      @.     @.     L'/L  @.     \\
@.      @.     @.    @VVV   @.     \\
0  @>>> L @>>> M @>>>  Q @>>> 0     \\
@.      @VVV    @|   @VVV     @.    \\
0  @>>> L' @>>> M @>>>  R @>>> 0     \\
@.      @.    @.   @VVV     @.    \\
@.      @.     @.     0     @.    \\
\endCD
$$
Note that $L'$, $M$ and $R$ satisfy the hypothesis of lemma \ref{bn2.6}, then
we can find (flat) deformations $M_y$ and $R_y$ parametrized by a
curve $Y$ such that for some $y_1\in Y$ we have that the support of the
corresponding sheaf $R_{y_1}$ is a smooth point of $C$. All sheaves and
maps can be deformed along.
To do this we define $Q_y=M_y/L$ (we have $L \inj L' \inj M_y$, thus
this quotient is well defined). The kernel of $Q_y \to R_y$ is $L'/L$. Then
$Q_y$ is a flat deformation (being the extension of a flat deformation
$R_y$ by a constant and hence flat deformation $L'/L$).
Then for each $y$ we have a commutative
diagram
$$
\CD
@.      @.     @.     0     @.     \\
@.      @.     @.    @VVV   @.     \\
@.      @.     @.     L'/L  @.     \\
@.      @.     @.    @VVV   @.     \\
0  @>>> L @>>> M_y @>>>  Q_y @>>> 0     \\
@.      @VVV    @|   @VVV     @.    \\
0  @>>> L' @>>> M_y @>>>  R_y @>>> 0     \\
@.      @.    @.   @VVV     @.    \\
@.      @.     @.     0     @.    \\
\endCD
$$
Observe that the length of the part of $Q_{y_1}$ supported in singular
points is $n-1$, so repeating this process we can deform $Q$ until its
support lies in the smooth part of $C$. This finishes the proof of the
proposition.

\end{proof}

\section{Proof of theorem I}
\label{bnProof of the main theorem}

In this section we will prove theorem I:

\begin{proof}
Nonemptyness follows from the fact that the Brill-Nother loci for
smooth curves is nonempty, and by upper semicontinuity of
$h^0(\cdot)$. By remark \ref{remark} we can assume $r>d-p_a$.
We will prove theorem I by applying lemma \ref{bn0.2}. 

We start
with a rank one torsion-free sheaf $A$ corresponding to a point in
$\overline W^r_d$, $d>0$, $r\geq 0$, with $\rho(r,d)>0$ (recall that
we are assuming $r>d-p_a$). We have $h^0(A)$, $h^1(A)>0$. As we
explained at the
beginning of the section \ref{bnGeneral case}, we call $A'$ its base point
free part. Then we take $B$ to be the base point free part of ${A'}^*$,
and finally define $A''$ to be $B^*$.

By lemma \ref{bn2.1}, $A''$ and ${A''}^*$ are rank one locally free
sheaves on $C$ generated by global
sections. Then by proposition \ref{bn1.5} we find a deformation
$\SA''$ of $A''$ parametrized by a some smooth irreducible curve $T$.

The support of $\SA''$ defines a family of curves $\SC$ parametrized
by the irreducible curve $T$. note that $\SC|_t$ is smooth for
$t\neq0$.

By the definition of $A'$ and $A''$ we have exact sequences
\begin{equation}
0 \to A' \to A \to Q \to 0
\label{short1}
\end{equation}
\begin{equation}
0 \to A' \to A'' \to \widetilde Q \to 0
\label{short2}
\end{equation}
with $h^0(A')=h^0(A)$ and $h^1(A'')=h^1(A')$ (lemma \ref{equalit}).
If we look at \ref{short2} we see that we are in the situation of
proposition \ref{bn2.3}, with $L=A'$, $M=A''$, $\SM=\SA''$, $P=T$.
Then we get a family $\SA'$ (parametrized by some connected but in
general not irreducible curve).
Now we use this family $\SA'$ and the sequence \ref{short1} to apply
\ref{bn2.5} with $L=A'$, $M=A$ and $\SL=\SA'$. We get a new family
$\SA$. We denote by $T'$ the curve parametrizing the family $\SA$.

This family satisfies all the hypothesis of lemma \ref{bn0.2}
(item (iv) is given by proposition \ref{bn2.2}), and
then theorem I is proved.

\end{proof}

\chapter{Irreducibility of the moduli space for $K3$ surfaces}
\label{k3}

In this chapter we will prove the following theorem:

\smallskip
\noindent\textbf{Theorem II.}
\textit{
With the notation of chapter \ref{Preliminaries}, if $L$ is a primitive
nonzero element of
$\Pic(S)$, and $H$ is an $(L,c_2)$-generic polarization, then
$\FM_H(L,c_2)$ is irreducible.}
\smallskip

Due to the fact that the moduli space is smooth, irreducibility is
equivalent to connectedness.

\bigskip
\textbf{Outline of the proof of theorem II}
\bigskip

First we will prove the theorem for the case in which Pic$(S)=\ZZ$
For $H$ to be $(L,c_2)$-generic we need $L$ to be an odd multiple of
a generator of Pic$(S)$, and tensoring the vector bundles with a line
bundle we can assume that $H=L$ is a generator of Pic$(S)$. After proving the
theorem for this case, in section \ref{General K3 surface} we show, 
by considering families of
surfaces, that if the result is true for Pic$(S)=\ZZ$, then it is also
true under the conditions of the theorem (this part is very similar to
an argument in \cite{G-H}). From now on we will assume
that Pic$(S)=\ZZ$ and that $H=L$ is the ample generator.

The proof is divided into two parts. In section \ref{Small second Chern 
class} we handle the case in which
$c_2 \leq \frac{1}{2} L^2 + 3$. First we see (proposition \ref{2.1}) that
the sheaves
satisfying this inequality are exactly those which are nonsplit
extensions of the form

$$ 0 \to \OOS \to V \to L \otimes I_Z \to 0,$$
with $l(Z)=c_2$. Then we study the set $X \subset \hilbcc$ for which there
exist nonsplit extensions like these above, and we see, 
using theorem I,
that it is connected (proposition \ref{2.2}). 
Finally we use this to
prove (proposition \ref{2.4}) the connectedness of $\modulil$ for 
$\dim \modulil>0$ (if the dimension is zero the result is known \cite{M}).

Note that for $c_2=\frac{1}{2} L^2 + 3$ we have $\dim \modulil =
L^2+6>0$, and then we can continue the proof by induction on $c_2$.

Let $C(n)$ be the set of irreducible components of $\FM(L,n)$. 
We construct a map

$$\Phi _n :C(n) \to C(n+1).$$
To define this map, take a sheaf $E$ in a component $A$ of $\FM(L,n)$.
Take a point $p \in S$ and a surjection $E \to \SO_p$. Let $F$ be the
kernel

$$0 \to F \to E \to \SO_p \to 0.$$
$F$ is clearly stable, and $c_2(F)=c_2(E)+1$. Now we define $\Phi _n
(A)$ to be the component in which $F$ lies. It is easy to see that
this is independent of all the choices made, so that $\Phi _n$ is well
defined.

Now we assume that $\modulil$ is irreducible for $c_2 < n$. We are
going to see that if every connected component of $\FM(L,n)$ has a
non-locally free sheaf $F$, $\Phi _{n-1}$ is surjective, and then by
induction $\FM(L,n)$ will be irreducible.

Let $B$ be a component of $\FM(L,n)$ with non-locally free sheaves. By
lemma \ref{3.3}, it has a non-locally free sheaf $F$ such that
$F^{\vee\vee} \in \FM(L,n-1)$. By smoothness of the moduli
space, $F^{\vee\vee}$ is in only one irreducible component.
Call this component $A$. By construction $\Phi_{n-1}(A)=B$.

In other words, we have seen that if $\FM(L,n-1)$ is irreducible, then
there is only one component $\FM_0$ of $\FM(L,n)$ that has sheaves
that are not locally free, and then to prove that the later has only
one component, it will be enough to check that every component has a
non-locally free sheaf.

We divide the possible values of $c_2$ in regions labeled by $n \geq
1$, with $c_2$ satisfying

$$
((n-1)^2 + (n-1) + \frac{1}{2}) L^2 + 3 < c_2 \leq (n^2 + n + \frac{1}{2}) 
L^2 + 3.
$$
If $V$ is locally free, we prove that then $V$ fits in a short exact
sequence

$$ 0 \to L^{\otimes -m} \to V \to L^{\otimes m+1} \otimes I_{Z_m} \to 0
$$
with $0 \leq m \leq n$ (proposition \ref{3.1}). We call it an
extension of type $m$. We will also say that $V$ is of type $m$.

Next (proposition \ref{3.3}) we show that the set of sheaves that are not
locally free has positive codimension, and then we prove (proposition  
\ref{3.4}) that the generic sheaf is a vector bundle of type $n$.

But this is not enough, and we need more information about the generic
vector bundle. Let $C$ be the set of vector
bundles $V$ such that for any exact sequence 

\begin{equation}
 0 \to L^{\otimes -n} \to V \to L^{\otimes n+1} \otimes I_{Z_n} \to 0,
\label{eq0.1}
\end{equation}
$L^{\otimes n+1} \otimes I_{Z_n}$ has no sections whose zero locus
is an irreducible reduced curve. In proposition \ref{3.7} we prove that
this set has positive codimension. The reason
to look at
this set is because it is precisely because of these sheaves that we
cannot apply the generalization of Fulton-Lazarsfeld's theorem to
proof that the set of type $n$ vector bundles is connected. But now we
know that we can ignore $C$, because it has positive codimension,
and then conclude that the generic vector bundle $V$ sits in an
extension like \ref{eq0.1} such that $L^{\otimes 2n+1} \otimes
I_{Z_n}$ has a section
whose zero locus is an irreducible reduced curve. In proposition \ref{3.6}
we prove that those vector bundles make a connected set. We will need
the induction hypothesis to prove this proposition.

\section{Preliminaries}
\label{secPreliminaries}

In this section we will prove some propositions that will be useful later.

\begin{lemma}
\label{1.1}
Let $S$ be a smooth surface and $C$ a smooth (not necessarily
complete)  curve. Let $p$ be a point
in the curve and $j:S \inj S \times C$ the corresponding injection.
Let $L$ be a line bundle on $S$ and $I_W$ an ideal sheaf on $S$
corresponding to a subscheme of dimension zero. Let $\SV$ be a family
of rank two sheaves on $S$, i.e. a sheaf on $S \times C$ flat over $C$. If we
have the following elementary transformation:

$$ 0 \to \SW \to \SV \to j_*(L \otimes I_W) \to 0 $$
then $\SW$ is a flat family of rank two sheaves on $S$, and furthermore 

$$ c_i(\SW_{p'}) = c_i(\SV_{p'})$$
for $i=1,2$ and $p'$ any point of $C$.
\end{lemma}

\begin{proof}
We calculate the Chern classes of $j_*(L \otimes I_W)$ by the
Grothendieck-Riemann-Roch theorem, and then the classes of $\SW$ by
Whitney's formula. The fact that $\SW$ is flat is proved in \cite{F}.
\end{proof}

Now we will apply this lemma to take limits of stable extensions. Let
$S$ be a smooth surface with $\Pic (S)=\ZZ$. Consider a family of
extensions parametrized by a curve $T$

$$ 0 \to L^{\otimes -n} \to V_t \to L^{\otimes n+1} \otimes I_Z \to 0,$$
where L is a generator of $\Pic (S)$, $t \in T$, and $Z$ is a subscheme
of dimension zero.
Assume that $V_t$ is stable for $t \not= 0$, where $0$ is some fixed
point of T, and unstable for $t = 0$. This defines a map $\varphi :
T-\{0\} \to \FM$ to the moduli space
of stable sheaves. By properness of $\FM$, this can be extended to a
map $\varphi : T \to \FM$, i.e. we can take the limit of the family as $t$
goes to $0$ and we obtain a stable sheaf corresponding to $\varphi (0)$.

\begin{proposition}
\label{1.2}
The stable sheaf $V'$ corresponding to $\varphi (0)$ is not locally free or can
be written as an extension

\begin{equation}
 0 \to L^{\otimes -m} \to V' \to L^{\otimes m+1} \otimes I_{Z'} \to 0
\label{eq1.1}
\end{equation}
with $m < n$.
\end{proposition}

\begin{proof}
$V_0$ is unstable, so we have

$$ 0 \to L^{\otimes a} \otimes I_W \to V_0 \to L^{\otimes 1-a} \otimes I_{W'}
\to 0$$
with $0 < a \leq n$. Consider the elementary transformation

$$ 0 \to \SW \to \SV \to j_*[L^{\otimes 1-a} \otimes I_{W'}] \to 0.$$
By lemma \ref{1.1} we have a new family $\SW$. By standard arguments the
member $W_0$ of the new family corresponding to $t=0$ can be written
as an extension

\begin{equation}
0 \to L^{\otimes 1-a} \otimes I_{W'} \to W_0 \to L^{\otimes a} \otimes I_W
\to 0. 
\label{eq1.2}
\end{equation}
Note that $0 \geq 1-a > -n$. If $W_0$ is not stable, repeat the
process: unstability gives an injective map $L^{\otimes a'}\otimes I_{W''}
\to W_0$, and by \ref{eq1.2} we have $a' < a$. Then in each step $1-a$ 
grows. We are going to see that eventually
we are going to get a stable sheaf. Assume we reach $1-a=0$ and $W_0$
is still unstable. The destabilizing sheaf has to be $L \otimes
I_{Z_d}$ with $l(Z_d)>l(W)$ and gives a short exact sequence

$$ 0 \to L \otimes I_{Z_d} \to W_0 \to I_{Z'_d} \to 0.$$
Note that $l(Z'_d)<l(W')$. Performing the corresponding
elementary transformation we
get a new family $\overline \SW$ and the sheaf corresponding to $0$
sits in an exact sequence

\begin{equation}
0 \to I_{Z'_d} \to \overline W_0 \to L\otimes I_{Z_d} \to 0.
\label{eq1.3}
\end{equation}
This sequence is like \ref{eq1.2} but with $0\leq l(Z'_d)<l(W')$. If we still
don't get a stable sheaf, repeat this. In each step $l(Z'_d)$
decreases, but this must stop because if $l(Z'_d)=0$, the sheaf
given by \ref{eq1.3} is stable, as the following lemma shows.

Now, once we have obtained a stable sheaf, if it is not locally free,
we are done. If it is locally free, then necessarily the subscheme
$W'$ is empty, and we get an extension
like \ref{eq1.1} as desired.
\end{proof}

\begin{lemma}
\label{1.3}
Let $V$ be a torsion free sheaf on a surface $S$ with $\Pic (S)=\ZZ$,
given by an extension

$$0 \to \OOS \to V \to L \otimes I_Z \to 0,$$
where $L$ is the effective generator of $\Pic (S)$. Then $V$ is
stable.
\end{lemma}

\begin{proof}
A destabilizing
subsheaf should be of the form $L^{\otimes m} \otimes I_W$, with $m>0$. By
standard arguments, it is enough to check stability with subsheaves
whose quotients are torsion free, so we can assume this.

The composition $L^{\otimes m} \otimes I_W \to V \to L \otimes I_Z$ is
nonzero, because otherwise it would factor through $\OOS$, but this is
impossible because $m>0$. Then $m=1$ and we have $I_W \inj I_Z$.
Furthermore, $l(W) > l(Z)$ because if $W=Z$, the sequence would split.

Then we have a sequence

$$ 0 \to L \otimes I_W \to V \to I_{W'} \to 0,$$
but we reach a contradiction because $c_2 = l(W) + l(W') > l(Z) +
l(W') = c_2 + l(W')$. Then there is no destabilizing subsheaf, and $V$
is stable.
\end{proof}

\begin{proposition}
\label{1.4}
Let $S$ be a smooth $K3$ surface with Picard group $\Pic (S)=\ZZ$. 
If $\dim \Ext ^1(L'
\otimes I_Z,L) \geq 2$, then there is a nonsplit extension 
\begin{equation}
0 \to L \to V \to L' \otimes I_Z \to 0 
\label{eq1.4}
\end{equation}
such that V is not locally free.
\end{proposition}

\begin{proof}
We have an exact sequence 

$$ 0 \to H^1(L\otimes(L')^{-1}) \to \Ext ^1(L' \otimes I_Z, L) \to H^0(\SO_Z).$$
If $L=L'$, then $H^1(L\otimes(L')^{-1})=0$ because $S$ is a K3 surface. If $L \neq
L'$, then due to the fact that $\Pic(S)=\ZZ$, applying Kodaira's
vanishing theorem we also have $H^1(L\otimes(L')^{-1})=0$. We have then an injection

$$ 0 \to \Ext ^1(L' \otimes I_Z,L) \stackrel{f}{\rightarrow} H^0(\SO_Z).$$
An extension corresponding to $\xi$ is locally free iff the section
$f(\xi)$ generates the sheaf $\SO_Z$, i.e., iff

$$ f(\xi) \notin W=\{ s \in H^0(\SO_Z) : 0=s \otimes k(p) \in
H^0(\SO_p) \text { for some } p \in \Supp(Z) \}.$$
$W$ is a union of codimension 1 linear subspaces, hence if $\dim 
\Ext ^1(L' \otimes I_Z, L) \geq 2$, then $\dim \text{im}(f) \cap
W > 0$, and we have a nonzero $\xi$ corresponding to an extension
\ref{eq1.4} with $V$ not locally free.
\end{proof} 

Usually we will apply the following corollary

\begin{corollary}
\label{1.5}
Let $S$ be a smooth $K3$ surface with $\Pic (S)=\ZZ$. If 
$$
\dim \Ext ^1 (L^{\otimes n+1}\otimes I_Z,L^{\otimes -n}) \geq 2,$$
 and there is a stable extension
$$0 \to L^{\otimes -n} \to V \to L^{\otimes n+1} \otimes I_Z \to 0,$$
then there is a sheaf $V'$, in the same irreducible component of
$\modulil$ as $V$, that is not locally free or sits in an extension

$$0 \to L^{\otimes -m} \to V \to L^{\otimes m+1} \otimes I_Z \to 0$$
for some $m<n$.
\end{corollary}

\begin{proof}
There is an open set
in $\PP(\Ext ^1 (L^{\otimes n+1} \otimes I_Z,L^{\otimes -n}))$ 
whose points correspond to
stable extensions, due to the openness of the stability condition. 
All these points get mapped to the same irreducible
component of $\modulil$. By proposition \ref{1.4}, there is an extension $V$
that is not locally free. If it is not stable, we can take a curve as
in proposition \ref{1.2}, and applying the proposition we get a family of
stable sheaves. All get mapped to the same component of $\modulil$,
and the sheaf corresponding to $t=0$ has the required properties.
\end{proof}

\section{Small second Chern class}
\label{Small second Chern class}

In this section we will consider the case in which $c_2 \leq \frac{1}{2}
 L^2 +3$. Recall that we are assuming that $S$ is a K3
surface with $\Pic (S)=\ZZ$. In this case we have the following
characterization of the stable torsion free sheaves.

\begin{proposition}
\label{2.1}
Let $V$ be a torsion free stable rank two sheaf with $c_1=L$, $c_2
\leq \frac{1}{2} L^2 + 3$. Then $V$ fits in an exact sequence

\begin{equation}
0 \to \OOS \to V \to L \otimes I_Z \to 0.
\label{eq2.1}
\end{equation}

Conversely, every nonsplit extension of $L \otimes I_Z$ by $\OOS$
is a torsion free stable sheaf.
\end{proposition}

\begin{proof}
Take $V$ stable. Using the Riemann-Roch theorem,

$$ h^0(V) + h^2(V) \geq \frac{L^2}{2} - c_2 + 4 \geq 1.$$
If $h^2(V)$ were different from zero, by Serre duality we would have
$\text {Hom} (V,\SO) \not= 0$, contradicting stability because this
would give a map $V \to \OOS$ with image $L^{\otimes -n} \otimes I_Z$ ($n
\geq 0$) and kernel $L^{\otimes n+1} \otimes I_{Z'}$.

Then $h^0(V) \not= 0$. Take a section of $V$. By stability, the
quotient of the section is torsion free, and we have an extension like
\ref{eq2.1}. The extension is not split because $V$ is stable.

The converse is lemma \ref{1.3}.
\end{proof}

Now that we know that all sheaves can be written as extensions of $L
\otimes I_Z$ by $\OOS$, the obvious strategy is to construct families
of extensions $\PP(\extlzo)$ for each $Z$ such that $\dim \extlzo \geq
1$. Ideally we would like to put all these families together in a
bigger family parametrized by a variety $M$. This $M$ would map to
$\modulil$ surjectively, so it would be enough to prove the
connectedness of $M$, and because $M$ maps to $\hilbcc$ with connected
fibers, it would be enough to prove that the set $X=\{ Z \in \hilbcc :
\dim \extlzo \geq 1\}$ (i.e., the image of the map $M \to \hilbcc$) is
connected.

Unfortunately we cannot construct $M$ because $\dim \extlzo$ is not
constant. We will use a somewhat more elaborate argument to bypass
this difficulty, but we will still use the connectivity of $X$, that
we prove in the following proposition.

\begin{proposition}
\label{2.2}
The set $X=\{ Z \in \hilbcc :\dim \extlzo \geq 1\}$ is connected.
\end{proposition}

\begin{proof}
By Serre duality and looking at the sequence

$$ 0 \to H^0(L \otimes I_Z) \to H^0(L) \to H^0(\SO_Z) \to H^1(L
\otimes I_Z) \to 0,$$
we have $\dim \extlzo \geq 1\ \iff h^0(L \otimes I_Z)\geq \frac{1}{2}
L^2 +3-c_2$. Now consider the following commutative diagram
$$
\CD
@.       @.            0                      @.  0                  @.\\
@.       @.          @VVV                       @VVV                 @.\\
0 @>>> \OOS @>>> L \otimes I_Z  @>>>  j_*(\omega _C \otimes I_Z) @>>> 0\\
@.       @|          @VVV                       @VVV                 @.\\
0 @>>> \OOS @>>>      L      @>>>       L| _C=j_*\omega _C @>>> 0\\
@.       @.           @VVV                       @VVV                @.\\
@.       @.          \SO_Z            @=         \SO_Z               @.\\
@.       @.           @VVV                      @VVV             @.\\
@.       @.           0                       @.  0                  @.\\
\endCD
$$
where $C \in \PP(H^0(L\otimes I_Z))$ (maybe $C$ is singular, but we know it is
irreducible and reduced because $\Pic (S)=\ZZ$ and $L$ is a generator of
the group), $j:C \inj S$ is the inclusion, and $\omega_C=L|_C$ is the
dualizing sheaf on $C$.

Using the top row we get $h^0(L \otimes I_Z)\geq \frac{1}{2}
L^2 +3-c_2 \iff h^0(\omega _C \otimes I_Z) \geq \frac{1}{2}
L^2 +2-c_2$. This condition can be restated in terms of Brill-Noether
sets $W^r_d$:

$$
\omega _C \otimes I_Z \in W^r_d
$$
where $r=\frac{1}{2}L^2 +1-c_2$, and $d=L^2 -c_2$.

By a theorem of Fulton and Lazarsfeld \cite {F-L}, the Brill-Noether 
set $W^r_d$
of a smooth curve is nonempty and connected if the expected dimension
$\rho (r,d) = g-(r+1)(g-d+r)$ is greater than zero. In the case of
an irreducible reduced curve lying on a $K3$ the generalized Jacobian can be
compactified, and the connectedness result is still true (theorem I). 
In our case we have

$$ \rho (r,d)=2c_2-\frac{L^2}{2}-3= \frac{{\dim \modulil}}{2} > 0,$$
(recall that for $\dim \modulil=0$ the irreducibility of the moduli space
is known by the work of Mukai \cite{M}) and we can apply the theorem. 
Now consider the variety 

$$N = \{ (Z,C): Z \subset C, \dim \extlzo \geq 1\} \subset \hilbcc 
\times \PP(H^0(L))$$
and the projections
$$
\CD
N @>p_2>> \PP(H^0(L))  \\
@Vp_1VV   @. \\
\hilbcc   @. \\
\endCD
$$

By theorem I, $p_2$ is surjective with connected
fibers. Then $N$ is connected, and also the image of $p_1$, that is
equal to $X$.
\end{proof}

For the following proposition we will need this lemma:

\begin{lemma}
\label{2.3}

Let $T$ be a smooth curve, $p$ a point in $T$ and $S$ a variety. Consider the
diagram
$$
\CD
S @>v>> S\times T \\
@VgVV  @VfVV \\
p @>u>> T \\
\endCD
$$
Then for every coherent torsion free sheaf $\SF$ on $S\times T$, there
exist a natural map
$$
(f_* \SF )(p) \to H^0(\SF _{S\times \{p\} }),
$$
where $\SF(p)=v^* \SF$. Furthermore, this map is injective.
\end{lemma}

\begin{proof}
The question is local in $T$, so we can assume that $T$ is affine,
$T=\Spec A$, and there is an element $x\in A$ such that $p$ is the
zero locus of $x$. We have
\begin{equation}
0 \to \SF \stackrel{\cdot x}{\to}  \SF \to \SF / x\cdot \SF \to 0,
\label{eq2.2}
\end{equation}
where the map on the left is multiplication by $f^*x$. On the other hand
we have the sequence
$$
0 \to \SO_{S \times T}(-f^*p) \to \SO_{S\times T} \to \SO_{S\times \{p\}
} \to 0. $$
Tensoring with $\SF$ is right exact, so we get an exact sequence
$$
\SF \otimes \SO_{S \times T}(-f^*p) \to \SF \to \SF_p \to 0.
$$
Note that the image of the left map is $x\cdot \SF$, and then we
conclude that $\SF / x\cdot \SF$ is isomorphic to $\SF_p$. Taking
cohomology in the sequence \ref{eq2.2} we get
$$
{H^0(\SF)}/(x\cdot H^0(\SF)) \inj H^0(\SF_p),
$$
but the first group is exactly $(f_* \SF)|_p$.
\end{proof}

Finally we can prove:

\begin{proposition}
\label{2.4}
The moduli space $\modulil$ of torsion free, rank two sheaves with
$c_2 \leq \frac{1}{2}L^2 + 3$ over a
K3 surface with $\Pic (S)=\ZZ$ is connected (hence irreducible, because
we know it is smooth).   
\end{proposition}

\begin{proof}
We have a stratification of $X$
$$
X=\bigcup _{r \geq 1} H_r,  H_r=\{ Z \in \hilbcc : \dim \extlzo
=r\}.
$$
On each stratum $H_r$ we can construct a projective bundle $M_r \to
H_r$ with fiber $\PP (\extlzo)$, because the dimension of the group is
constant. Each point of $M_r$ corresponds to an extension (up to weak
isomorphism of extensions). We have then morphisms $M_r \to  \FM
(L,c_2)$ with fiber $\PP(H^0(V))$ over $V$ (see proposition \ref{3.8}). We have
$h^0(V)=\frac{1}{2}L^2 +3-c_2 + h^1(L \otimes I_Z)$, and a
corresponding stratification of $\modulil$
$$
\modulil=\bigcup _{r \geq 1} \FM _r ,\ \ \ 
\FM _r = \{V \in \FM (L,c_2) : h^0(V) = \frac{1}{2}L^2 +3-c_2
+ r \}
$$
(the reason for this dependence on $r$ in the definition is that $h^0(V)
= \frac{1}{2}L^2 +3-c_2+h^1(L\otimes I_Z)$. The condition is
equivalent to requiring that if $V \in
\modulil$ and $V$ is an extension of $L\otimes I_Z$ by $\OOS$, then
$h^1(L\otimes I_Z)=r$. The previous formula proves that $r$ only
depends on $V$).
 
Note that $M_r$ can be thought also as a projective bundle over $\FM
_r$ with fiber $\PP(H^0(V))$.

To prove that $\modulil$ is connected, we will show that for any two
sheaves $V_a \in \FM _a$, $V_b \in \FM _b$, we can construct a family
$\SV$ of stable sheaves on a connected parameter space, with $\SV |_0
= V_a$, $\SV |_1 = V_b$

Due to the fact that $X$ is connected, it is enough to prove this for
$V_a$, $V_b$ given by extensions
$$
\CD
0 \to \OOS @>s_a>> V_a \to L \otimes I_{Z_a} \to 0 \\
0 \to \OOS @>s_b>> V_b \to L \otimes I_{Z_b} \to 0
\endCD
$$
such that $Z_a \in H_a$, $Z_b \in H_b$ and $Z_a \in \overline H_b$,
the closure of $H_b$.

Take a curve $f:T \to \hilbcc$ with $f(0)=Z_a$, $f(1)=Z_b$ and
$\im (T-\{0\}) \in H_b$. $T$ doesn't need to be complete. By shrinking
$T$ to a smaller open set, we can assume that there is a lift
$f$ to a map $f:T\setminus\{0\} \to M_b$. This gives a family of sheaves $\widetilde
\SV$ and sections $s_t \in H^0(\SV|_t)$ parametrized by $T\setminus\{0\}$.

We want to extend this to a
family $\SV$ of sheaves and sections parametrized by $T$, in such a way 
that the
cokernel of $s_0:\OOS \to \SV |_0$ is $L \otimes I_{Z_a}$. Maybe $\SV |_0$
won't be isomorphic to $V_0$, but at least both are extensions of the
same sheaf $L \otimes I_{Z_a}$ by $\OOS$, and then they are in the
same connected component, and this is enough.

This family gives a morphism $T-\{0\}\to \modulil$ that extends to a
unique $T \to \modulil$ by properness (see section \ref{secPreliminaries}, 
just before
proposition \ref{1.2}). With this we have already
extended the family $\widetilde \SV$ to a family $\SV$ parametrized by
$T$, and we only need to extend the sections $s_t$. Shrinking $T$ to a
smaller neighborhood of $\{0\}$ if
necessary, we can assume that $\pi _2{}_* \SV$ is trivial ($\pi _1$
and $\pi _2$ are the first and second projections of $S \times T$).
Then the sections $s_t$ fit together to give $\SO_{S \times (T-\{0\})}
\to \widetilde \SV$, i.e. an element
$\tilde s \in ({\pi_2} _* \SV)(T-\{0\})$. This can be extended to some
$s\in ({\pi_2} _* \SV)(T)$ that is nonzero on the fiber of $t=0$.

Using the previous lemma, we have an injection $({\pi_2} _* \SV)|_{t=0}
\to H^0(\SV |_0)$ that gives a nonzero section $s_0$ of $\SV |_0$. Now
we only have to check that the cokernel of this section is $L \otimes
I_{Z_a}$. We have a short exact sequence over $S \times T$
$$
0 \to \SO_{S \times T} \stackrel{s}{\to} \SV \to \SQ \to 0.
$$

Then $\SQ$ is torsion free, flat over $T$, and then it is of the form
$$\SQ = {\pi _1}^* (L) \otimes I_\SZ \otimes {\pi _2}^* (L')$$ 
for some line bundle $L'$ over $T$ and some subscheme $\SZ$ of $S \times T$ 
flat over $T$. This subscheme gives a morphism $g:T \to \hilbcc$. 
By construction
$f(t)=g(t)$ for $t \not= 0$, and by properness the equality also holds
for $t=0$, then $Z_0 = Z$ as desired.
\end{proof}

\section{Large second Chern class}
\label{Large second Chern class}

In this section we will handle the case in which $c_2$ is large.

\begin{proposition}
\label{3.1}
Assume that $c_2$ satisfies
$$
((n-1)^2 + (n-1) + \frac{1}{2}) L^2 + 3 < c_2 \leq (n^2 + n + \frac{1}{2}) 
L^2 + 3.
$$
If $V$ is locally
free, then $V$ fits in a short exact sequence
$$
0 \to L^{\otimes -m} \to V \to L^{\otimes m+1} \otimes I_{Z_m} \to 0,
l(Z_m) = c_2 + (m+1) m L^2 ,
$$
with $0 \leq m \leq n$. We will call such an exact sequence an
extension of type $m$.
\end{proposition}

\begin{proof}
For any sheaf $V$,
$h^2(V\otimes L^{\otimes n}) = \dim \Hom(V, L^{\otimes -n}) = 0$ by stability,
and then using the Riemann-Roch theorem we have 
\begin{equation}
h^0(V\otimes L^{\otimes n}) \geq \frac{L^2}{2} - c_2 + 4 +(n + n^2) L^2 \geq 1,
\label{eq3.1}
\end{equation}
so that we have an inclusion $L^{\otimes -n} \inj V$. If $V$ is locally
free, this will give an exact sequence on type $m$, with $m \leq n$.
\end{proof}

\begin{proposition}
\label{3.2}
If a component $\FM'$ has sheaves of type $m$, with $0 \leq m \leq
n-1$, then it has sheaves that are not locally free.
\end{proposition}

\begin{proof}
Choose $m$ such that there is no $V$ of type $m'$ for $m' < m$. Let
$V$ be of type $m$. By Serre duality
$$
\dim \Ext ^1(L^{\otimes m+1} \otimes I_{Z_m}, L^{\otimes -m}) = h^1(L^{\otimes 2m+1}
\otimes I_{Z_m}),
$$
and this is greater than 2. By proposition \ref{1.5}, $\FM'$ has a sheaf
that is not locally free or is of type $m'<m$, but the later cannot
happen because of the choice of m.
\end{proof}  

\begin{lemma}
\label{3.3}
Let $X$ be an irreducible component of $\modulil$. Let $X_0$ be the
subset corresponding to non-locally free sheaves. If $X_0$ is not
empty, then it has codimension one. Furthermore, there is a dense
subset $Y$ of $X_0$ such that for any $F \in Y$, we have
$$
c_2(F^{\vee\vee})=c_2(F)+1.$$
\end{lemma}

\begin{proof}
By \cite {O1}, prop. 7.1.3, we know that 
\begin{equation}
\codim(X_0,X)\leq 1.
\label{bound}
\end{equation}
On the other hand, let $F\in X_0$. It fits in an exact sequence
$$
0 \to F \to F^{\vee\vee} \to \SQ \to 0
$$
where $\SQ$ is an Artinian sheaf
with length $l=H^0(\SQ)=c_2(F^{\vee\vee})-c_2(F)$. We use this to
bound the dimension of $X_0$ by a parameter count. First we choose a
locally free sheaf $E\in \FM(L,c_2-l)$. These requires
$4(c_2-l)-L^2-6$ parameters. Now we have to choose a quotient to a
sheaf of length $l$ concentrated on a subset of dimension zero. These
quotients are parametrized by the Grothendieck Quot scheme $\Quot
(E,l)$, whose dimension is $3l$ (this follows from \cite {L1}
Appendix, where it is proved that $\Quot ^0(E,l)$, the Quot scheme
corresponding to quotients supported in $l$ distinct points, is dense
in $\Quot(E,l)$).

Define the following stratification on $X_0$:
$$
X_0=\bigcup_{l\geq 1} X_0^l \text{,   }X_0^l=\{F:
c_2(F^{\vee\vee})-c_2(F)=l\}
$$
Then we have
$$
\dim (X_0^l) \leq 4(c_2-l)-L^2-6 + 3l
$$
and together with the previous bound \ref{bound} we obtain that
$Y=X_0^1$ is dense
and $\codim(X_0,X)=1$.
\end{proof}

\begin{proposition}
\label{3.4}
If $c_2$ satisfies the inequalities of the hypothesis of proposition
\ref{3.1}, then there is an open dense set on $\modulil$ that
corresponds to sheaves of type n.
\end{proposition}

\begin{proof}
We will prove this by showing that the codimension of sheaves of type
$m \leq n-1$ is greater than zero.

We will divide the proof into two cases:

\textbf{Case 1.}
\textit{Extensions of type $m$ with $Z_m$ such that $h^0(L^{\otimes 2m+1} \otimes
I_{Z_m}) = 0$.}

We have then $h^1(L^{\otimes 2m+1} \otimes I_{Z_m}) = c_2 -(m^2 +m+\frac{1}{2}) 
L^2 -2$. The dimension of the family $M$ of extensions of this kind is
bounded (via Serre duality) by
$$
\dim M \leq 2 l(Z_m) + h^1(L^{\otimes 2m+1} \otimes I_{Z_m}) -1 = 3 c_2
+(m^2 +m-\frac{1}{2}) L^2 -3.
$$
There is a map $\pi:M \to \modulil$ with fiber over each $V$ equal to
$\PP(H^0(V\otimes L^{\otimes n}))$, and we can give a bound to its dimension (see proof
of proposition \ref{3.1}).
$$
\dim \PP(H^0(V\otimes L^{\otimes n})) \geq (n^2 + n + \frac{1}{2}) L^2 - c_2 + 3
$$
Then the
dimension of the image of $\pi$ is bounded by 

\begin{eqnarray*} 
\dim (\im \pi)  & \leq & \dim (M) - \min \dim \PP(H^0(V\otimes L^{\otimes n})) \\
                & \leq & 4 c_2 - L^2 -6 - 2nL^2
\end{eqnarray*}

And then
$$
\codim (\im \pi) \geq 2nL^2 >0
$$

\textbf{Case 2.}
\textit{Extensions of type $m$ with $Z_m$ such that 
$h^0(L^{\otimes 2m+1} \otimes I_{Z_m}) \geq 1$.}

Let $M$ be the family of stable extensions with $Z_m$ satisfying this
inequality. Let $l = l(Z_m) = c_2 + (m+1)mL^2$.
The subscheme $Z_m$ is in a curve $C$ defined as the zeroes of a section
of $L^{\otimes 2m+1}$, i.e. $Z_m
\in \Hilb ^l  (C)$. Although this curve will be reducible and
not reduced in general, the fact that $C$ is in a smooth surface
allows us to prove:

\begin{lemma}
\label{3.5} 
$\dim \Hilb ^d (C) = d$
\end{lemma}

\begin{proof}
We have a natural stratification of $\Hilb ^d (C)$ given by the number
of points in the support of a subscheme
$$
\Hilb ^d (C) = \bigcup _{1 \leq r \leq d} \Hilb ^d _r (C),
$$
where $\Hilb ^d _r (C) = \{ Z : \# \Supp Z = r\}$. We have natural maps
giving the support of a subscheme:
$$
z^C _r : \Hilb ^d _r(C) \to \Sym ^r(C).
$$

In the same way we define maps $z^S _r$, when we consider subschemes
of the surface $S$. Clearly, if $x \in \Sym ^r (C) \subset \Sym ^r
(S)$, then $(z^C _r)^{-1} (x) \subset (z^S _r)^{-1} (x)$. Taking this
into account, and using a result of Iarrobino about zero dimensional
subschemes of a smooth surface \cite {I}:

\begin{eqnarray*}
\dim \Hilb ^d _r (C) &\leq & \dim \Sym ^r (C) + \maxdim (z^C _r)^{-1}
(x) \\
                     &\leq & \dim \Sym ^r (C) + \maxdim (z^S _r)^{-1} (x)\\
                     &=    & r + d - r = d
\end{eqnarray*}

This gives $\dim \Hilb ^d (C) \leq d$, and the opposite direction is
trivial. This finishes the proof the lemma.
\end{proof}

Now we are going to bound the dimension of the set 
$$
H' = \{ Z \in \Hilb^l (S): h^0(L^{\otimes 2m+1} \otimes I_Z) \geq 1\}.
$$
Consider the diagram
$$
\CD
Y = \{ (Z,C) \in \Hilb^l \times \PP(H^0(L^{\otimes 2m+1}) : Z \in C\}
@>p_2>>
\PP(H^0(L^{\otimes 2m+1})) \\
@Vp_1VV      @. \\
\Hilb ^l (S) @. \\
\endCD
$$
We have $H' = \im (p_1)$ and the fiber of $p_2$ is $\Hilb ^l (C)$. $p_2$
is clearly surjective and the
fiber of $p_1$ is $\PP (H^0(L^{\otimes 2m+1} \otimes I_Z))$.  Then
$$
\dim H' = \dim \PP (H^0(L^{\otimes 2m+1})) + \dim \Hilb ^l (C) - \dim \PP
(H^0(L^{\otimes 2m+1} \otimes I_Z)) = 2l-1.$$

Again we have a map $\pi:M \to \modulil$ with fiber
$\PP (H^0(V\otimes L^{\otimes m}))$, and then
$$
\codim (\im \pi) \geq (2m+1) L^2 +1 > 0.
$$
\end{proof}

As a corollary to this proposition we learn that to prove
connectedness of $\modulil$ it is enough to prove that all type $n$
sheaves are in one component.

\begin{proposition}
\label{3.6}
All stable extensions of $L^{\otimes n+1} \otimes I_Z$ by $L^{\otimes -n}$ such that
$L^{\otimes 2n+1} \otimes I_Z$ has
a section corresponding to an integral curve, are in one component.
\end{proposition}

\begin{proof}
Define the sets
\begin{eqnarray*}
\widetilde X_r & = & \{Z \in \hilbn : \dim \extn = r \text { and} \\
               &  &L^{\otimes 2n+1} \otimes I_Z \text { has a section 
corresponding to an integral curve} \} \\
 & & \\
X_r &= &\{ Z \in \widetilde X_r : \text { there is a stable extension of
} L^{\otimes n+1} \otimes I_Z \text { by } L^{\otimes -n} \} \\
\widetilde M_r &=& \{ \text {extensions of } L^{\otimes n+1} \otimes I_Z 
\text {by } L^{\otimes -n} \text { with } Z \in X_r \}\\
 & & \\
M_r &=& \{ m \in M_r : m \text { corresponds to a stable extension} \} \\
 & & \\
N_r &=& \{(Z,C) \in \hilbn \times \PP(H^0(L^{\otimes 2n+1})) : Z \subset C  \\
 & & \text { and } \dim \extn = r \} \\
 & & \\
U &=& \{ C \in \PP(H^0(L^{\otimes 2m+1})) : C \text { is irreducible and reduced} \}
\end{eqnarray*}

We construct $\widetilde M_r$ as parameter spaces of universal families of
extensions by standard techniques. These techniques require that the
dimension of the $\Ext ^1$ group is constant on the whole family. This
why we have to introduce the subscript $r$ and break
everything into pieces according to the dimension of the group $\Ext
^1$. We also consider the unions
$$
\widetilde M = \bigcup _{r \geq 1} \widetilde M_r \text {,  }
M = \bigcup _{r \geq 1} M_r \text {,  }
X = \bigcup _{r \geq 1} X_r \text {, } \ldots
$$

Note that $X$, being a subset of $\hilbn$, has a natural scheme 
structure. This is also true for  $N \subset C
\times \PP(H^0(L^{\otimes 2n+1}))\ $. On the other hand, for $M_r$ 
there is no natural way of
``putting them together'', so we take just the disjoint union. 
We have the following maps
$$
\CD
\widetilde M     @<<<   M   @>>>  \modulil \\
@VVV                          @VVV         @.        \\
\widetilde X     @<<<   X        @.        \\
@A{p_1}AA                     @.           @.        \\
N     @>{p_2}>>          \PP(H^0(L^{\otimes 2n+1}))      @.    \\
\endCD
$$ 

By construction $\widetilde X = p_1 {p_2}^{-1}(U)$.
Now we are going to prove that the fibers of $p_2$ over $U$ are
nonempty and connected. For each point in $N$ we have a commutative
diagram
$$
\CD
@.       @.            0                      @.  0                  @.\\
@.       @.          @VVV                       @VVV                 @.\\
0 @>>> \OOS @>>> L^{\otimes 2n+1} \otimes I_Z @>>> 
j_*(\omega _C \otimes I_Z) @>>> 0\\
@.       @|          @VVV                       @VVV                 @.\\
0 @>>> \OOS @>>>  L^{\otimes 2n+1}    @>>> 
L^{\otimes 2n+1}|_C=j_*\omega _C @>>> 0\\
@.       @.           @VVV                       @VVV                @.\\
@.       @.          \SO_Z            @=         \SO_Z               @.\\
@.       @.           @VVV                      @VVV             @.\\
@.       @.           0                       @.  0                  @.\\
\endCD
$$

We argue in the same way as in proposition \ref{2.2}. Here we have
$$
r=(n^2+n+\frac{1}{2})L^2 +1 -c_2 \text{, }
d=(3n^2+3n+1)L^2-c_2
$$
$$ \rho (r,d)=2c_2-\frac{L^2}{2}-3=\frac{{\dim \modulil}}{2} > 0.$$
But now we don't know if $p_2$ is surjective with connected fibers,
because theorem I only applies
for irreducible reduced curves. This is the reason why we introduce
the open set $U$. For the fibers on $U$ we can apply the theorem, and
we conclude that $p_2$ is surjective over $U$ with connected fibers.
Then $\widetilde X=p_1 {p_2}^{-1}(U)$ is connected.

\textbf{Case 1.}
\textit{$p_1 {p_2}^{-1} (U) = X$.}

If $X_1 = X$, then we can construct a (connected, because $X$ is
connected) family $M_1$ parametrizing all sheaves with the
required properties, and we are done.

Now, if $X_1 \not= X$, then there are extensions with $r \geq 2$. By
corollary \ref{1.5}, $M_r$ with $r \geq 2$ is mapped to $\FM _0$, the
irreducible component that has sheaves that are not locally free. There is 
only one irreducible component with this property, because by induction
hypothesis the moduli space when the second Chern class is smaller 
than $c_2$ is 
irreducible (see the outline of the proof). Now
we have to show that all the connected components of $M_1$ also go to
this component $\FM _0$.

The connectivity of $X=\widetilde X$ and the fact that $\dim \extn$ is
upper semicontinuous allows us to take a curve $f:T \to X$ with
$f(T-\{0\})$ in any given connected component of $X_1$ and $f(0) \in
X_r$, with $r \geq 2$.

Lift $f$ to a map $f:T-\{0\} \to M_1$. Note that $M_1$ won't be in
general a projective bundle because we have removed the points
corresponding to unstable extensions, but these make a closed set, and
(maybe after restricting $T$ to a smaller open set) we can construct
the lift without hitting this set.

$M_1$ maps to $\modulil$, and then we have a map $\phi : T-\{0\} \to
\modulil$. As in the proof of proposition \ref{2.4}, this gives us a
family of stable sheaves and sections parametrized by $T-\{0\}$ that
we can extend to a family parametrized by $T$. Now there are two
possibilities:

If $\phi(0)$ is of type $n$, then we have a family of extensions
$$
0 \to L^{\otimes -n} \to V_t \to \OOS(L^{\otimes n+1} \otimes I_{Z_t}) \to 0 $$
and a corresponding map $\psi : T \to \hilbn$, $t \mapsto Z_t$. By
construction $\psi(t) = f(t)$ for $t \not= 0$, and by properness also
for $t=0$. The extension corresponding to $t=0$ has to be in $M_r$
with $r \geq 2$, and then $M_1$ is also mapped to $\FM _0$.

On the other hand, if $\phi(0)$ is not a vector bundle of type $n$,
then it is either of type $m$ for $m<n$ or it is not locally free.
In either case, we conclude that $M_1$ is also mapped to $\FM _0$

\textbf{Case 2.}
\textit{$p_1 {p_2}^{-1} (U) \not= X$.}

Again, $M_r$, $r \geq 2$, gets mapped to $\FM _0$.

No connected component of $X_1$ can be closed, because by
connectedness of $\widetilde X$ and upper semicontinuity of $\dim
\extn$, we would have $\widetilde X = X_1$, and then $X = \widetilde
X$, contrary to the hypothesis.

Now we can prove that every connected component of $M_1$ gets mapped
to $\FM _0$. Take the corresponding connected component of $X_1$. Take
a curve $f:T \to \widetilde X$, with $f(T-\{0\})$ in the given
connected component of $X_1$, and $f(0) \notin X_1$. As in the
previous case, lift $f$ to a map $f:T-\{0\} \to M_1$, and now the proof
finishes like the end of case 1.
\end{proof}

\begin{proposition}
\label{3.7}
The set of sheaves $V$ of type $n$ such that for any extension 
$$
0 \to L^{\otimes -n} \to V \to L^{\otimes n+1} \otimes I_Z \to 0,$$
$L^{\otimes 2n+1} \otimes I_Z$ has no section whose zero locus is an
integral curve, has positive codimension. 
\end{proposition}

\begin{proof}
Define $\widetilde P=\{ Z \in \hilbn : L^{\otimes 2n+1} \otimes I_Z$ has no
sections whose
zero locus is an irreducible reduced curve $\}$.
For each point of $\widetilde P$ we have a family of extension of type $n$
given by the projectivization of the corresponding $\Ext ^1$ group.
Writing $\widetilde P = \cup \widetilde P_r$, with $r$ equal to the
dimension of the group, we can
construct a family of extensions $\widetilde M_r^P$ for each $r$. As is
the previous proposition, let $P_r \subset \widetilde P_r$ be the
subset that has stable extensions.
We have a natural map $\pi_1:M_r^P \to \modulil$, where $M_r^P$ is the
subset of $\widetilde M_r^P$ corresponding to stable sheaves.

\begin{lemma}
\label{3.8}
The fiber of $\pi_1$ over $V \in \modulil$ is 
$\PP (H^0(V\otimes L^{\otimes n}))$.
\end{lemma}

\begin{proof}
The fiber consists of all extensions giving the same $V$. Now, given a
point in $\PP (H^0(V\otimes L^{\otimes n}))$, we have an injection 
$f: L^{\otimes -n} \inj V$ (up
to scalar). $V$
is locally free and of type $n$, then the quotient is torsion free
and we get an element $Z_f$ of $\hilbn$:

$$0 \to L^{\otimes -n} \to V \to L^{\otimes n+1} \otimes I_{Z_f} \to 0.$$
This defines a map from $\PP (H^0(V\otimes L^{\otimes n}))$ to the fiber of $\pi$. It is
clearly surjective. Now we will check that it is also injective.

If $Z_f = Z_{f'}$, then $f$ and $f'$ have to differ at most by scalar
multiplication, because all nonsplit extensions of $L^{\otimes n+1} \otimes
I_{Z_f} = L^{\otimes n+1} \otimes I_{Z_{f'}}$ by $L^{\otimes -n}$ that give the same $V$
are weak isomorphic, so we get a diagram: 
$$
\CD
0 @>>> L^{\otimes -n} @>f>> V @>>> L^{\otimes n+1} \otimes I_{Z_f} @>>> 0 \\
@.  @V{\alpha}VV @V{\isom}VV @V{\beta}VV @. \\
0 @>>> L^{\otimes -n} @>f'>> V @>>> L^{\otimes n+1} \otimes I_{Z_{f'}} @>>> 0 \\
\endCD
$$
where $\alpha$ is multiplication by scalar.
\end{proof}

In $\PP (H^0(L^{\otimes 2n+1}))$ we have a subvariety $Y$ corresponding to
reducible curves. This subvariety is the image of the natural map
$$\bigcup _{0<a<2n+1} \PP(H^0(L^{\otimes a})) \times 
\PP(H^0(L^{\otimes 2n+1-a})) \to
\PP(H^0(L^{\otimes 2n+1}))$$
We define the set

\begin{eqnarray*}
\widetilde N_r^P & = &\{ (Z,C) \in X_r \times \PP(H^0(L^{\otimes 2n+1})): \\
                 &   &Z \subset C, \dim (p_1)^{-1}(Z)=(n^2+n+
\frac{1}{2})L^2-c_2 +1+r \}
\end{eqnarray*}

and the maps
$$
\CD
\widetilde N_r^P @>p_2>> \PP(H^0(L^{\otimes 2n+1}))  \\
@Vp_1VV   @. \\
\hilbn   @. \\
@A{\pi_2}AA   @. \\
M_r^P @>{\pi_1}>> \modulil \\
\endCD
$$

By construction we have $P_r \subset \widetilde P_r \subset
p_2 p_1^{-1} (Y)$. Then
$$
\dim P_r \leq \dim p_2((p_1)^{-1} (Y)) = \dim Y + \dim\fiber (p_2) -
\dim\fiber (p_1),
$$
where $\dim Y$ is the maximum of the dimensions of its irreducible
components. Finally
$$\codim (\im \pi_1) = \dim \modulil - \dim P_r -\dim\fiber \pi_2 +
\dim\fiber \pi_1,$$
and putting everything together we have $\codim (\im \pi_1) >
(2n-a)(a-1)L^2$ for every $0<a<n+1$, and then $\codim (\im \pi_1) > 0$.
\end{proof}

\section{General K3 surface (proof of theorem II)}
\label{General K3 surface}

In this section we finally prove theorem II by showing that if the
result is true for a surface $S$ with $\Pic(S)=\ZZ$, then it also
holds under the hypothesis of theorem II. The idea is to deform the
given surface to a generic surface with $\Pic(S)=\ZZ$. We also deform
the moduli space, and then the irreducibility of the moduli space for
the deformed surface will imply the irreducibility for the surface we
started with. This is very similar to an argument in
\cite{G-H}. 

Because we are going to vary the surface, in this
section we will denote the moduli space of semistable sheaves with
$\FM_H(S,L,c_2)$, where $S$ is the surface on which the sheaves are
defined.

\smallskip
\begin{proof2} \textit{of theorem II}
Recall that we have a surface $S$ with a $(L,c_2)-$generic polarization $H$.
By 2.1.1 in \cite{G-H}, there is a connected family of surfaces $\calS$
parametrized by a curve $T$ and a line bundle $\SL$ on $\calS$ such that 
$(\calS_0,\SL_0) = (S,L)$ and $\Pic(\calS_t)=\SL_t \cdot \ZZ$
for $t\neq 0$. By proposition 2.3 in \cite{G-H}, there is a connected smooth
proper family $\SZ \to T$ such that $\SZ_0 \isom \FM_H(S,\SL_0,c_2)$
(note that the polarization is $H$ and not $\SL_0$)
and  $\SZ_t \isom \FM_{\SL_t}(\calS_t,\SL_t,c_2)$ for $t \neq 0$.

By propositions \ref{2.4} and \ref{3.4} we know that $\SZ_t$ is irreducible
for $t \neq 0$, and then by an argument parallel to lemma \ref{bn0.2}, we 
obtain that $\SZ_0$ is connected, but $\SZ_0$ is smooth
(because $H$ is generic), and then this implies that $\SZ_0$ is 
irreducible.
\end{proof2}

\chapter{Irreducibility of the moduli space for del Pezzo surfaces}
\label{dp}

In this chapter we will consider the case in which $S$ is a del Pezzo
surface. We will prove the following theorem (see chapter
\ref{Preliminaries} for the notation).

\smallskip
\noindent\textbf{Theorem III.}
\textit{
Let $S$ be a del Pezzo surface. Let $\FM_L(S,c_1,c_2)$ be the moduli
space of rank 2, Gieseker L-semistable torsion free sheaves with Chern
classes $c_1$, $c_2$, with $L$ a $(c_1,c_2)$-generic polarization.
Then $\FM_L(S,c_1,c_2)$ is either empty or irreducible.}
\smallskip

As we explained in chapter \ref{Preliminaries}, this is already known
for $\PP^2$ and $\PP^1 \times \PP^1$, so we will assume from now on
that $S$ is a surface isomorphic to $\PP^2$ blown up at 
most at 8 points in general position. We denote the blow up map 
$\pi:S \to \PP^2$. We will denote by $H$ the effective generator of
$\Pic(\PP^2)$. We also denote by $H$ the pullback $\pi^*(H)$.
$\FM_L(X,c_1,c_2)$ will denote the moduli space of $L$-semistable
torsion free rank two sheaves on $X$ with Chern classes $c_1$, $c_2$.
In the case $X=\PP^2$ we will drop $L$ from the notation, because
there is only one possible polarization.

This is a particular case of the conjecture of Friedman and Qin
\cite{F-Q} that states that the moduli
space of stable vector bundles on a rational surface with effective
anticanonical line bundle is either empty or irreducible (for any
choice of polarization and Chern classes).

\begin{proposition}
\label{dp1}
$\FM_{L_0}(S,c_1,c_2)$ is either irreducible or empty, where
$L_0$ is a polarization of $S$ that lies in a chamber whose closure
contains $H$.
\end{proposition}

\begin{proof}
By \cite{B}, $\FM_{L_0}(S,bH+a_1E_1+\cdots+a_nE_n,c_2)$ is
birational to a $(\PP^1)^m$ bundle over $\FM(\PP^2,bH,c_2)$ (where $m$ 
is the number of $a_i$'s that are odd), but
it is well known that this moduli space is either irreducible or
empty.
\end{proof}

To prove this statement for any generic polarization, we will
need some lemmas about the following system of equations on integer
numbers:

$$
\left.
\begin{array}{rl}
a_1^2+ \cdots +a_8^2 &= x+b^2\\
-a_1- \cdots -a_8&=x-2+3b
\end{array}
\right \}
\ \ \ 
\begin{array}{c}
(\dagger)
\end{array}
$$
where $x$ is some given (integer) number, and $b$,$a_1,\ldots,a_8$
are the unknowns.

\begin{lemma}
\label{dp2}
If $x\geq3$ then any integer solution of $(\dagger)$ with 
$b>0$ has $b\leq2$.
\end{lemma}

\begin{proof}
This is obtained by an elementary argument.
We can interpret geometrically these equations as the intersection of
a one-sheeted hyperboloid and a plane. We will look first at real
solutions. Rewrite $(\dagger)$ as
$$
\left.
\begin{array}{rl}
b &= \frac {1}{3}(-a_1 - \cdots -a_8 +2-x) \\
a_1^2+ \cdots +a_8^2 &= x+\frac {1}{9}(-a_1 - \cdots -a_8 +2-x)^2
\end{array}
\right \}
$$
The first equation defines a function, and the second equation is a
constrain. Using the method of Lagrange multipliers we obtain that the
maximum and minimum values of $b$ are at points of the form $a_i=t$
for some $t$. Then $(\dagger)$ become
$$
\left .
\begin{array}{rl}
8t^2 & =b^2+x \\
-8t-3b & =-2+x
\end{array}
\right \}
$$
Looking at the real solutions of these equations, we find
\begin{equation}
\label{bpm}
b^\pm = \frac {-6(x-2)\pm \sqrt{36(x-2)^2-4((x-2)^2-8x)}}{2}.
\end{equation}
For $x\geq 3$, $b^-$ is always negative. 
If we want to have solutions with $b>0$
we need $b^+>0$. Using \ref{bpm}, this is equivalent to
$(x-2)^2<8x$, and this implies $x<(12+\sqrt{128})/2<12$. This bound,
together with the hypothesis $3\leq x$ and \ref{bpm} implies
$b^+<3$, but if we are only interested in integer solutions this gives
$b\leq 2$.
\end{proof}

\begin{lemma}
\label{dp3}
If $1\leq b \leq 2$, then the integer solutions of $(\dagger)$ (up to
permutation of $a_i$) are:
\begin{eqnarray*}
b=1,\ a_1=\cdots=a_{x+1}=-1,\ a_{x+2}=\cdots=a_8=0 & \\
b=2,\ a_1=\cdots=a_{x+4}=-1,\ a_{x+5}=\cdots=a_8=0 & \! .
\end{eqnarray*}
\end{lemma}

\begin{proof}
In both cases ($b=1$ or $b=2$), the right hand sides of the two
equations are equal, and then, substracting the equations we have
$$
\sum a_i^2+a_i=0.
$$
But for $a_i$ integer we have $a_i^2+a_i\geq0$, and then $a_i$ must be
equal to -1 or 0.

Now, looking at the equations we see that the number of nonzero
$a_i$'s is given by $x+b^2$, and we obtain the result.
\end{proof}

\begin{theorem}
\label{dp4}
$\FM_L(S,c_1,c_2)$, for any generic polarization $L$, is either empty
or irreducible.
\end{theorem}

\begin{proof}
We will denote by $L_0$ a polarization lying in a chamber whose
closure contains $H$.

If $\FM_L(S,c_1,c_2)$ has more than one irreducible component, then
there must be a wall between $L$ and $L_0$ that created the extra
component. Recall that a wall $W^\zeta$ is a hyperplane in the ample
cone perpendicular to a class $\zeta$ with $\zeta \equiv c_1$ (mod 2),
and $c_1^2-4c_2\leq\zeta^2<0$. By \cite{F-Q}, if $L_1\cdot\zeta>0
>L_2\cdot\zeta$, then the sheaves that are $L_1$-unstable and
$L_2$-stable make an irreducible family of dimension $N_\zeta+2l_\zeta$, where
$N_\zeta=h^1(\zeta)+l_\zeta-1$ and $l_\zeta=(4c_2-c_1^2+\zeta^2)/4$.
In the case of a rational surface we have
$$
h^1(\zeta)=\frac{\zeta\cdot K_S}{2}-\frac{\zeta^2}{2}-1.
$$
The wall creates a new component if $N_\zeta+2l_\zeta$ is equal to the
dimension of the moduli space, in our case $4c_2-c_1^2-3$. For a
rational surface we have $N_\zeta+N_{-\zeta}+2l_\zeta=4c_2-c_1^2-4$,
and then this condition is equivalent to $N_{-\zeta}=-1$. Denoting
$\zeta=bH+a_1E_1+\cdots+a_nE_n$ and $x=4c_2-c_1^2$ we get the system
of equations ($\dagger$) ($n\leq8$ so without loss of generality we can study
the equations with $n=8$). Furthermore, this wall will create a new 
component if
$L_0\cdot\zeta<0<L\cdot\zeta$. By the definition of chamber, the last
equality is equivalent to $0<H\cdot\zeta$, and this
translates to $b>0$.

We will prove the proposition by showing that for given $(S,c_1,c_2)$
there is at most one such wall in the ample cone and that in this
case $\FM_{L_0}(S,c_1,c_2)$ is empty, so that $\FM_L(S,c_1,c_2)$ is
always either empty or only has one irreducible component.

If the moduli space is not empty, its dimension should be
nonnegative, and this translates to $x\geq3$. Then by lemmas \ref{dp2}
and \ref{dp3} we know all the solutions of ($\dagger$), i.e. all the
walls creating components.

By tensoring with a line bundle (and relabeling the exceptional
curves), we can assume that $c_1$ is either $H+E_1+\cdots+E_m$ or
$E_1+\cdots+E_m$ for some $m\leq8$. Now we will use the condition $c_1
\equiv \zeta$ (mod 2). 

In the first case this implies that the only
possible solution for ($\dagger$) is $\zeta=H-E_1-\cdots-E_m$, and $m=x+1$.
Then $4c_2-c_1^2=x$ implies $c_2=0$. By \cite{B},
$\FM_{L_0}(S,H+E_1+\cdots+E_m,0)$ is birational to a $(\PP^1)^m$ bundle over
$\FM(\PP^2,H,0)$, but it is well known that this moduli space is
empty, then the same is true for $\FM_{L_0}(S,H+E_1+\cdots+E_m,0)$.

In the second case ($c_1=E_1+\cdots+E_m$), the condition $c_1 \equiv
\zeta$ (mod 2) implies that the only possible solution is
$\zeta=2H-E_1-\cdots-E_m$ with $m=x+4$. Then $c_2=-1$ and then
$\FM_{L_0}(S,c_1,c_2)$ is empty (because its expected dimension is
negative).
\end{proof}
 
This techniques can also be used to study the irreducibility of the
moduli space of stable torsion free sheaves on $\PP^2$
with more than 8 blown up points. We obtain equations similar to
$(\dagger)$, but with more variables. Unfortunately now we don't have
a bound on $b$ like the one given by lemma \ref{dp2}, and then it
becomes more difficult to classify the solutions. This is still work
in progress and it will appear elsewhere.

\end{document}